\def\VERSION{5.7.2026}
\def\WHO{Tom} 
\def\users{us}    
\def\users{world} 
\documentclass[12pt]{article} 

\usepackage{color}
\usepackage{amsmath}
\usepackage{amsthm}
\usepackage{amssymb}
\usepackage{epsfig}
\usepackage{psfrag}
\usepackage{graphicx}
\usepackage{textcomp}
\usepackage{bm}
\usepackage{pgf}
\numberwithin{equation}{section}
\usepackage{upgreek} 
\newtheorem{theorem}{Theorem}[section]

\newtheorem{definition}[theorem]{Definition}

\newtheorem{proposition}[theorem]{Proposition}

\newtheorem{remark}[theorem]{Remark}

\newcommand{\ITEM}[2]{\parbox[t]{.055\textwidth}{{\rm #1}}\hfill\parbox[t]{.945\textwidth}{#2}\vspace*{.8mm}}

\usepackage{mathrsfs,cite}
\usepackage{ifthen}
\usepackage{ulem}
\usepackage[colorlinks=true, pdfstartview=FitV, linkcolor=blue, 
            citecolor=blue, urlcolor=blue]{hyperref}
\usepackage{cancel}\ifthenelse{\equal{\users}{world}}
{
\newcommand{\REM}[1]{}

	\newcommand{\DELETE}[1]{}

        \newcommand{\COMMENT}[1]{}
        \newcommand{\TCOMMENT}[1]{}
    \newcommand{\MARGINOTE}[1]{}
}	
{
\usepackage[notcite,notref,color]{showkeys}

\definecolor{brown}{rgb}{0.6,0.2,0.2}
\newcommand{\REM}[1]{\marginpar{\bfseries\tiny{\color{blue}#1}}}

 \newcommand{\COMMENT}[1]{{\color{blue}\uuline{#1}\color{black}}} 
 \newcommand{\DELETE}[1]{{\color{brown}\cancel{#1}\color{black}}}

 \newcommand{\TCOMMENT}[1]{{\color{blue}{ #1}}}
\newcommand{\MARGINOTE}[1]{\marginpar{\color{red}\tiny\texttt{#1}}}
\usepackage{fancyhdr}
\pagestyle{fancy}
\headheight=28pt\headwidth=17cm
\definecolor{gray}{gray}{0.5}
\newcount\hour \newcount\minute
\hour=\time
\divide \hour by 60
\minute=\time
\loop \ifnum \minute > 59 \advance \minute by -60 \repeat

\rhead{\color{gray}time discretization Eulerian\\
by T.Roub\'\i\v cek}
\chead{}
\lhead{Version \VERSION, file: \jobname.tex, \underline{\Large\color{red}\WHO's working on} \\
compiled:
\number\day.\number\month.\number\year\ at
\the\hour:\ifnum\minute<10 0\fi\the\minute\ h }
}

\newcommand{\R}{\mathbb{R}}

\newcommand{\bbI}{\mathbb{I}}

\newcommand{\bbD}{\mathbb{D}}

\newcommand\DT[1]{\mathchoice
                 {{\buildrel{\hspace*{.1em}\text{\LARGE.}}\over{#1}}}
                 {{\buildrel{\hspace*{.1em}\text{\LARGE.}}\over{#1}}}
                 {{\buildrel{\hspace*{.1em}\text{\Large.}}\over{#1}}}
                 {{\buildrel{\hspace*{.1em}\text{\large.}}\over{#1}}}}

\newcommand{\lineunder}[2]{\LU{\begin{array}[t]{c}\underbrace{#1}\vspace*{.5em}\end{array}}{\mbox{\footnotesize\rm #2}}}
\newcommand{\linesunder}[3]{\LSU{\begin{array}[t]{c}\underbrace{#1}\vspace*{.5em}\end{array}}{\mbox{\footnotesize\rm #2}}{\mbox{\footnotesize\rm#3}}}
\newcommand{\LU}[2]{\begin{array}[t]{c}#1\vspace*{-1em}\\_{#2}\end{array}}
\newcommand{\LSU}[3]{\begin{array}[t]{c}#1\vspace*{-1em}\\_{#2}\vspace*{-.5em}\\_{#3}\end{array}}

\newcommand{\threelinesunder}[4]{\threeLSU{\begin{array}[t]{c}\underbrace{#1}\vspace*{.5em}\end{array}}{\mbox{\footnotesize\rm #2}}{\mbox{\footnotesize\rm #3}}{\mbox{\footnotesize\rm #4}}}
\newcommand{\threeLSU}[4]{\begin{array}[t]{c}#1\vspace*{-1em}\\_{#2}\vspace*{-.5em}\\_{#3}\vspace*{-.5em}\\_{#4}\end{array}}

\renewcommand{\d}{{\rm d}}

\newcommand{\divS}{\mathrm{div}_{\scriptscriptstyle\textrm{\hspace*{-.1em}S}}^{}}
\newcommand{\eq}[1]{(\ref{#1})}
\newcommand{\Cdot}{\hspace{-.1em}\cdot\hspace{-.1em}}
\newcommand{\Colon}{\hspace{-.15em}:\hspace{-.15em}}

\def\vv{{\bm v}}
\def\pp{{\bm p}}

\def\zz{{\bm z}}
\def\AA{{\bm A}}

\def\xx{{\bm x}}
\def\yy{{\bm y}}
\def\nn{{\bm n}}

\newcommand{\DD}{\bm D}
\newcommand{\TT}{\bm T}
\newcommand{\XX}{\bm X}
\newcommand{\Se}{\bm S}
\newcommand{\strain}{{\boldsymbol\varepsilon}}
\newcommand{\GRAVITY}{\bm g}

\def\vvk{\vv_\etau^k}

\def\overlinemuetau{\hspace*{.2em}\overline{\hspace*{-.2em}\mu}_\etau^{}}
\def\overlinealphaetau{\hspace*{.2em}\overline{\hspace*{-.2em}\alpha}_\etau^{}}
\def\underlinealphaetau{\hspace*{.1em}\underline{\hspace*{-.1em}\alpha}_\etau^{}}
\def\overlineEetau{\hspace*{.2em}\overline{\hspace*{-.2em}\bm E}_\etau^{}}

\def\overlinevvtau{\hspace*{.15em}\overline{\hspace*{-.15em}\vv}_{\etau}^{}}
\def\overlinemutau{\hspace*{.15em}\overline{\hspace*{-.15em}\mu}_{\etau}^{}}
\def\ovetau{\hspace*{.15em}\overline{\hspace*{-.15em}\vv}_{\etau}^{}}
\def\opetau{\hspace*{.15em}\overline{\hspace*{-.15em}\pp}_{\etau}^{}}

\def\overlinexitau{\hspace*{.15em}\overline{\hspace*{-.15em}\bm\xi}_{\etau}^{}}
\def\underlinexitau{\hspace*{.0em}\underline{\hspace*{-.0em}\bm\xi}_{\etau}^{}}
\newcommand{\Fezero}{\FF_{\hspace*{-.1em}\mathrm e,0}^{}}
\newcommand{\Fekk}{\FF_{\hspace*{-.1em}\mathrm e,\tau}^{k-1}}
\newcommand{\Fem}{\FF_{\hspace*{-.1em}\mathrm e,\tau}^{m}}
\newcommand{\Femm}{\FF_{\hspace*{-.1em}\mathrm e,\tau}^{m-1}}

\def\FF{{\bm F}}
\def\Fe{{\FF_{\rm e}}}
\def\Fee{{\FF_{\rm e}}}
\def\Fp{{\FF_{\rm p}}}
\def\Fpp{{\FF_{\rm p}}}
\def\MM{{\bm M}}
\newcommand{\zetad}{\zeta_{\rm dm}}
\newcommand{\zetap}{\zeta_{\rm vp}}
\def\Lp{{\bm L}}

\def\HYPER{\nu}
\def\eetau{\tau}

\newcommand{\Frac}[2]{\mathchoice{\text{\small$\frac{#1}{#2}$}}
                                 {\text{\large$\frac{#1}{#2}$}}
                                 {\text{\large$\frac{#1}{#2}$}}
                                 {\text{\large$\frac{#1}{#2}$}}}
\newcommand{\barOmega}{\,\overline{\!\varOmega}}
\newcommand{\nablaS}{\nabla_{\scriptscriptstyle\textrm{\hspace*{-.3em}S}}^{}}

\newcommand{\Nabla}{\nabla}
\newcommand{\Rsym}{\mathbb R^{3\times3}_{\rm sym}}
\newcommand{\Rdev}{\mathbb R^{3\times3}_{\rm dev}}

\def\Vdots{\!\mbox{\setlength{\unitlength}{1em}
\begin{picture}(0,0)
\put(-.07,0){.}
\put(-.07,.3){.}
\put(-.07,.6){.}
\end{picture}
}
}

\newcommand{\wt}[1]{\mathchoice{\hspace*{-.09em}\text{\large$\hspace*{.09em}\tilde{\text{\normalsize$#1$}}\hspace*{.05em}$}\hspace*{-.05em}}
{\hspace*{-.09em}\text{\large$\hspace*{.09em}\tilde{\text{\normalsize$#1$}}\hspace*{.05em}$}\hspace*{-.05em}}
{\text{\normalsize$\hspace*{.08em}\tilde{\text{\scriptsize$#1$}}\hspace*{.06em}$}}
{\text{\small$\tilde{\text{\tiny$#1$}}$}}}

\newcommand\DELETEDELETE[1]{}

\newcommand\pdt[1]{\frac{\partial{#1}}{\partial t}}

\begin{document}

\def\EPS{\varepsilon}
\def\DELTA{\delta}
\def\etau{{\tau}}
\def\EEps{{}}

\def\TTtauk{\TT^k_{\!\etau}}
\def\overlineDetau{\hspace*{.2em}\overline{\hspace*{-.2em}\bm D}_\etau^{}}
\def\overlineTetau{\hspace*{.2em}\overline{\hspace*{-.1em}\bm T}_\etau^{}}
\def\overlineLpetau{\hspace*{.2em}\overline{\hspace*{-.3em}\bm\varPi\hspace*{-.1em}}_\etau^{}}
\def\overlineMetau{\hspace*{.4em}\overline{\hspace*{-.4em}\bm M\hspace*{-.1em}}_\etau^{}}
\def\overlineFetau{\hspace*{.2em}\overline{\hspace*{-.2em}\bm F}_{\!{\rm e},\etau}^{}}
\def\underlineFetau{\hspace*{.0em}\underline{\bm F\hspace*{-.2em}}_{\,{\rm e},\etau}^{}}
\def\overlineHetau{\hspace*{.2em}\overline{\hspace*{-.2em}\bm H}_\etau^{}}
\def\overlineEetau{\hspace*{.2em}\overline{\hspace*{-.2em}\bm F}_{\!{\rm e},\tau}^{}}
\def\overlinevvtau{\hspace*{.1em}\overline{\hspace*{-.1em}\vv}_{\etau}^{}}
\def\overlineppetau{\hspace*{.15em}\overline{\hspace*{-.15em}\pp}_\etau^{}}
\def\overlinerhoetau{\hspace*{.15em}\overline{\hspace*{-.15em}\varrho}_\etau^{}}
\def\oretau{\hspace*{.15em}\overline{\hspace*{-.15em}\varrho}_\etau^{}}
\def\osetau{\hspace*{.15em}\overline{\hspace*{-.15em}\sigma}_\etau^{}}

\newcommand{\Fek}{\FF_{\hspace*{-.1em}\mathrm e,\tau}^k}
\newcommand{\Feek}{\Fee_{\hspace*{-.1em}\mathrm e,\tau}^k}
\newcommand{\overlineFe}{\hspace*{.2em}\overline{\hspace*{-.2em}\FF}_{\hspace*{-.2em}\mathrm e,\tau}^{}}

\allowdisplaybreaks

\noindent{\LARGE\bf Staggered time discretization 
\\[.2em]in finitely-strained heterogeneous 
\\[.2em]visco-elastodynamics with damage \\[.2em]
or diffusion in the Eulerian frame.
}

\bigskip\bigskip

\noindent{\large\sc Tom\'{a}\v{s} Roub\'\i\v{c}ek}\footnote{{\tt ORCID}: 0000-0002-0651-5959.}
\\
{\it Mathematical Institute, Charles University, \\Sokolovsk\'a 83,
CZ--186~75~Praha~8,  Czech Republic
}\\and\\
{\it Institute of Thermomechanics, Czech Academy of Sciences,\\Dolej\v skova~5,
CZ--182~08 Praha 8, Czech Republic
}

\bigskip\bigskip

\begin{center}\begin{minipage}[t]{14.5cm}

{\small

\noindent{\bfseries Abstract.}
The semi-implicit (partly decoupled, also called staggered or fraction-step)
time discretization is applied to compressible nonlinear
dynamical models of viscoelastic solids in the Eulerian description, i.e.\ in
the actual deforming configuration, formulated fully in terms of rates. The
Kelvin-Voigt rheology and also, in the deviatoric part, the Jeffreys rheology
are considered. The numerical stability and, considering the Stokes-type
viscosity multipolar of the 2nd-grade, also convergence towards weak solutions
are proved in three-dimensional situations, exploiting the convexity of the
kinetic energy when written in terms of linear momentum.
No (poly)convexity of the stored energy is required and some enhancements
(specifically towards damage and diffusion models) are briefly outlined, too.

\medskip


\medskip


} 
\end{minipage}
\end{center}

\bigskip

\section{Introduction}

The {\it visco-elastodynamics} models in continuum mechanics {\it at finite}
(also called {\it large}) {\it strains} lead to strongly nonlinear systems of
evolution partial differential equations. In solid mechanics, the Lagrangian
approach (using a referential configuration) is most commonly used, but
sometimes the {\it Eulerian approach} (using the actual deforming
configuration and the return mapping to facilitate heterogeneous media)
has specific advantages. In particular, sometimes any referential configuration
is not justified (e.g., in geophysical models over long time scales of
millions of years) and interactions with external spatial fields (such
as gravitational or electromagnetic)
or fluid-solid interactions is more straightforward in the Eulerian frame.
Moreover, the Stokes viscosity is simpler in the Eulerian frame than
in the Lagrangian frame where frame-indifference requires special
nonlinearities and special analytical effort.

Using the  Eulerian  approach, the goal of this article is to modify
the fully implicit time discretization devised and analyzed for the homogeneous
polyconvex stored energies in \cite[Sect.3]{Roub25TDVE}. More specifically,
we use a suitable truncation of the stored energy as devised already in
\cite{Roub24TVSE,RouSte23VESS,RouTom23IFST}, where it was used for
a ``semi-Galerkin'' approximation however, and then devise a partly
{\it decoupled time discrete scheme}.  In addition, another ``nonlocal''
truncation of the stress is devised in Section~\ref{sec-diff-ext}. 
This substantially widens the
applicability of this time-discretization method to various general problems,
including those with stored energies depending on the space (i.e.,
{\it heterogeneous} problems) and on some internal variables (as {\it damage} or
a {\it diffusant concentration}) in a general nonconvex way.  The importance
of such problems is reflected by a vast literature about numerical treatment
of compressible media, mainly of an engineering character but also mathematical,
see in particular the monographs \cite{SNPeOw08CMPT,FeFeSt03MCMC,FLMS21NACF,FeKaPo16MTCV,LewSch98FEMS,SimHug98CI}.

As the visco-inertial part of the model is the same as
\cite{Roub25TDVE}, we will benefit from referring to several technicalities
to it. Yet, in contrast to \cite{Roub25TDVE}, we consider a rather general
(not necessarily quadratic) viscoplastic dissipation
potential  and the general (not necessarily polyconvex) stored energy.

The core of the (partially) decoupled time discretization is the coupled
system for mass density and velocity but with a truncated conservative
Cauchy stress discretized explicitly from the previous time step. Neglecting
the elastic 
response, this core would reduce to the compressible viscoelastic fluid system
for which such implicit time discretization has sometimes been used, mostly
merely computationally while for an  analytically rigorous treatment we
refer to
\cite{GaMaNo19EEIM,FHMN17EENM,FeKaPo16MTCV,FLMS21NACF,Kar13CFEM,Zato12ASCN}.
This documents the extreme analytical difficulty, which often requires strong
qualifications of the data or of solutions as strict positivity of mass
density in \cite{FLMS21NACF} or two space  dimension  in
\cite{Zato12ASCN} etc.
Here the elastic shear response and coupling with viscoplastic
Jeffreys type rheology, spatially nonhomogeneous (i.e.\ heterogeneous)
situations, and the mentioned coupling with other internal variables
would further emphasized such analytical difficulties.
For this analytical reason, involving some higher-order gradients
(here in the Eulerian model in the dissipative
part) facilitates the rigorous analysis and, in addition,
yields some additional parameter for a possible tuning of dispersion of the
velocity of elastic wave propagation, cf.\ \cite{Roub24SGTL}.
These higher (here 2nd-order) gradients lead to the concept of (here 2nd-grade)
{\it nonsimple media}, which often occurs in literature since the works by
R.A.\,Toupin \cite{Toup62EMCS} and R.D.\,Mindlin \cite{Mind64MSLE}. In the
dissipative part as used in this paper, it was also developed
in particular in \cite{BeBlNe92PBMV,FriGur06TBBC,NeNoSi89GSIC,NecRuz92GSIV}
as {\it multipolar} fluids.

We will use the standard notation concerning the Lebesgue and the Sobolev
spaces of functions on the Lipschitz bounded domain $\varOmega\subset\R^3$,
namely $L^p(\varOmega;\R^n)$ for Lebesgue measurable $\R^n$-valued functions
$\varOmega\to\R^n$ whose Euclidean norm is integrable with $p$-power
with $1\le p<\infty$, while $p=\infty$ means essentially bounded functions. 
For $k$ integer, $W^{k,p}(\varOmega;\R^n)$ stands for functions from
$L^p(\varOmega;\R^n)$ whose all derivatives up to the order $k$ have their
Euclidean norm integrable with $p$-power.
We have the embedding $W^{1,p}(\varOmega)\subset L^\infty(\varOmega)$ for $p>3$.
The star-notation $(\cdot)^*$ will be used for the dual space. In particular
$L^p(\varOmega)^*=L^{p'}\!(\varOmega)$ with the conjugate exponent
$p'=p/(p{-}1)$ if $1<p<+\infty$ or $p'=\infty$ if $p=1$. Moreover,
for a Banach space $X$ and for $I=[0,T]$, we will use the notation
$L^p(I;X)$ for the Bochner space of Bochner measurable functions $I\to X$
whose norm is in $L^p(I)$.
The space of continuous functions on the closure $\barOmega$ of
$\varOmega$ will be denoted by $C(\barOmega)$.

For readers' convenience, let us summarize the basic notation used in what
follows: 
\begin{center}
\fbox{
\begin{minipage}[t]{17em}\small\smallskip
$\yy$ deformation,\\
$\vv$ velocity,\\
$\varrho$ mass density,\\
$\sigma=1/\varrho$ sparsity,\\
$\pp=\varrho\vv$ the linear momentum,\\
$\FF=\nabla\yy$ deformation gradient,\\ 
$\Fe$ elastic distortion,\\
$\Fp$ inelastic distortion,\\
$\alpha$ damage or diffusant content,\\
$\bm\xi$ return mapping,\\
$\GRAVITY$ gravity acceleration,\\
$\HYPER$ the hyper-viscosity coefficient,\\
$I=[0,T]$ a time interval, $T>0$,\\
$\R_{\rm sym}^{3\times3}$ set of symmetric matrices,\\
``$\,\Vdots\,$'' scalar products 3rd-order tensors,\\ 
$(\cdot)'$ (partial) derivative of a  mapping,\\[.1em]
$(\cdot)\!\DT{^{\,}}$ convective time derivative,
\end{minipage}
\begin{minipage}[t]{21em}\small\smallskip
$\TT$ Cauchy stress,\\
$\DD$ dissipative stress,\\
$\MM$ Mandel's stress,\\
$\Lp$ inelastic distortion rate,\\
$\mu$ chemical potential,\\ 
$m$ mobility of diffusant (diffusivity),\\ 
tr$(\cdot)$ trace of a matrix,\\
dev$(\cdot)$ deviatoric part of a matrix,\\
$\varphi$ stored energy,\\
$\zetap$ viscoplastic dissipation potential,\\
$\zetad$ dissipation potential for damage,\\
$\bbI$ the identity matrix,\\
$\R_{\rm dev}^{3\times3}=\{A\in\R_{\rm sym}^{3\times 3};\ {\rm tr}A=0\}$,\\
${\rm GL}_3^+=\{A\in\R^{3\times3};\ \det A>0\}$,\\
``$\:\Cdot\:$'', ``$\:\Colon\:$'' scalar products of vectors or matrices,\\
$(\cdot)^*$ a dual space or a convex conjugate,\\
$\tau>0$ a time step for discretization.
\end{minipage}
}\end{center}

\vspace{-.8em}

\begin{center}
{\small\sl Table\,1.\ }
{\small
Summary of the basic notation used. 
}
\end{center}

The plan of this paper is the following: in Sect.\,\ref{sec-system}
the visco-elastodynamic model is formulated as a boundary-value problem
for a system of semi- and quasi-linear equations. Then, in
Sect.\,\ref{sec-discretization}, we then devise a truncated de-coupled
(called also staggered or fractional-step) time discretization, proved
its stability and convergence (in terms of subsequences) to weak solutions
of the previously formulated model. Finally, some extensions
of the model facilitated by the staggered discretization strategy
are briefly outlined in Sect.\,\ref{sec-multiphysics}, illustrated
specifically on a damage model and on a diffusion in poroelastic
media. 

\section{The visco-elastodynamic system and its energetics}\label{sec-system}

We demonstrate the essence of the staggered time-discretization method on
the visco-elasto\-dynamics in {\it Kelvin-Voigt rheology} in the volumetric
part and a {\it Jeffreys} (also known as {\it anti-Zener}) {\it rheology} in
the deviatoric (isochoric) part, which is a fairly general model that
allows for isochoric (viscoplastic) creep.

\subsection{The kinematics in brief}

In the large-strain continuum mechanics, the fundamental geometrical concept is a 
{\it deformation}\index{deformation} $\yy:\varOmega\to\R^3$ as a mapping from a
reference configuration $\varOmega\subset\R^3$ into the physical space $\R^3$.
The inverse motion $\bm\xi=\yy^{-1}:\yy(\varOmega)\to\varOmega$, if it
exists, is called a {\it return } {\it mapping}. We will denote the reference
(Lagrangian) and the actual
(Eulerian) point coordinates by $\XX$ and $\xx$, respectively. The further
basic geometrical object is the (referential)
{\it deformation gradient}\index{deformation!gradient}
$\FF_\text{\!\sc r}^{}(\XX)=\Nabla_{\XX}^{}\yy$.

If evolving in time, $\xx=\yy(t,\XX)$ is sometimes called a ``motion''.
The important quantity is the (referential) velocity
$\vv_\text{\sc r}^{}=\frac{\d}{\d t}\yy(t,\XX)$ with $\d/\d t$ the derivative
with respect to time of a time dependent function. When composed with the
return mapping $\bm\xi$, we obtain the Eulerian representations 
$\FF(t,\xx)=\FF_\text{\!\sc r}^{}(t,\bm\xi(\xx))$ and 
$\vv(t,\xx)=\vv_\text{\!\sc r}^{}(t,\bm\xi(\xx))$.
The Eulerian velocity $\vv$ is employed in the convective time derivative
$(\bm\cdot)\!\DT{^{}}=\pdt{}(\bm\cdot)+(\vv\Cdot\nabla)(\bm\cdot)$
with $\nabla$ taken with respect to actual
coordinates, to be used for scalars and, component-wise, for vectors or tensors.

The following transport equation by mere advection holds for the return mapping:
\begin{align}\label{ultimate-}
\DT{\bm\xi}=\bm0\ \ \ \text{ with the initial condition }\
\bm\xi|_{t=0}^{}=\bm\xi_0^{}\,.
\end{align}

Then the velocity gradient
$\Nabla\vv=\nabla_{\!\XX}^{}\vv\nabla_{\!\xx}^{}\XX=\DT\FF\FF^{-1}$,
where we used the chain-rule calculus  and
$\FF^{-1}=(\nabla_{\!\XX}^{}\xx)^{-1}=\nabla_{\!\xx}^{}\XX$. 
This gives the {\it transport-and-evolution equation for the
deformation gradient} as
\begin{align}\nonumber\\[-2.2em]
\DT\FF=(\nabla\vv)\FF\,.
  \label{ultimate}\end{align}
The deformation $\yy$ itself
does not need to be explicitly involved in the model when formulated in rates;
in fact, $\yy$ does not even need to exist but, if it does exist,
then $\FF=\nabla\yy$; cf.\ \cite[Remark~7]{Roub22VELS}.

To implement Maxwellian-type visco-elastic rheology, one should introduce
a certain decomposition of the total deformation gradient. Most commonly,
the Kr\"oner-Lee-Liu \cite{Kron60AKVE,LeeLiu67FSEP} {\it multiplicative
decomposition} is used, i.e.\
\begin{align}\label{KLL}
\FF=\Fe\Fp\ \ \ \text{ with $\ \Fe\ $ and $\ \Fp\ $ the elastic and the
inelastic distortions, }
\end{align}
respectively. The interpretation of $\Fp$ is a transformation of the reference
configuration into an intermediate stress-free configuration, and then the
{\it elastic distortion} $\Fe$ transforms this
intermediate configuration into the current actual configuration.
The stored energy as introduced below then depends on $\Fe$, if no
isotropic hardening is considered.

Applying the material derivative to \eq{KLL} and using \eq{ultimate}, we obtain
$(\Nabla\vv)\FF=\DT\FF=(\Fe\Fp)\DT{\!\!^{}}\!\!=\DT\Fe\Fp+\Fe\DT\Fp$
and, multiplying it by $\FF^{-1}=\Fpp^{-1}\Fee^{-1}$, eventually we obtain
the  decomposition $\Nabla\vv
=\DT\Fe\Fee^{-1}+\Fe\Lp\Fee^{-1}$ if we define the {\it inelastic distortion rate}
as $\Lp=\DT\Fp\,\Fpp^{-1}$, which is the most commonly used modelling concept.
From this, we obtain 
\begin{align}
\DT\Fe=(\Nabla\vv)\Fe+\Fe\Lp\,.
\label{dafa-formula}\end{align}
Noteworthy, we can eliminate $\Fp$ and formulate the evolution
merely in terms of $\vv$ and $\Fe$ if the stored energy will not
involved $\Fp$. Moreover, it is natural to consider the fluidic type
rheology only in the shear part while the volumetric part should be of a
solid type, for which $\det\Fp$ should be constrained, most commonly
equal 1, i.e.\ $\Fp$ is isochoric. This non-linear  constraint 
$\det\Fp=1$ is easily ``translated'' into a linear constraint ${\rm tr}\Lp=0$
thanks to the evolution $\DT\Fp=\Lp\Fp$ which then ensures
$(\det\Fp)\!\DT{^{}}=(\det\Fp){\rm tr}\Lp=0$, while assuming (rather
implicitly) that the initial
condition is $\det\Fp|_{t=0}=1$.

Without going into the details of the above standard concepts, we refer
to many renowned monographs, e.g.\ to \cite{GuFrAn10MTC,Mart19PCM,Rubi21CMEF}.

\subsection{The visco-elastic rheology and the dynamics}

The basic ingredients for the Maxwellian-type rheology,
which combines the (hyper)elastic and the viscous responses
in series, are the stored energy $\varphi=\varphi(X,F)$
and viscoplastic dissipation potential $\zetap=\zetap(X,L)$ depending 
on the position $X\in\varOmega$, the elastic distortion $F\in\R^{3\times3}$,
and the inelastic distortion rate $L\in\R^{3\times3}$. Moreover, 
from the analytical reasons, we consider also another viscosity
of the Stokes type parallel to the hyperbolic-type Maxwellian
model, which eventually make the visco-elastic rheology (called
Jeffreys' or anti-Zener's) of the parabolic type. 

The stored energy is considered in the physical unit J/m$^3$=Pa, with m$^3$
referring to the actual volume, not the referential volume. This gives
the (symmetric) {\it Cauchy stress} tensor $\varphi_F'(X,F)F^\top\!+\varphi(X,F)$.
Let us only mention that other alternatives, such as energy per referential
volume or per mass, would lead to the Cauchy stress of the form
$\varphi_F'(X,F)F^\top\!/\det F$ or $\rho(X)\varphi_F'(X,F)F^\top\!$,
respectively; here $\rho$ denotes the mass density at the referential
coordinate. The {\it Maxwellian}-type {\it viscosity} governed
by the potential $\zetap$ is driven by the deviatoric part
of the {\it Mandel stress} tensor $M$, i.e.\ dev$\,M=M-({\rm  tr}\,M)\bbI/3$
with $M=F^\top\varphi_F'(X,F)$. The {\it mass conservation} is 
reflected by the continuity equation $\pdt{}\varrho={\rm div}(\varrho\,\vv)$
for the actual mass density $\varrho$. 

Rigorously, we notationally have
distinguished the placeholder $F$ for values of the field $\Fe$, and
similarly $L$ versus $\Lp$, or $X$ versus  $\XX=\bm\xi(\xx)$.

Merging the mentioned mass conservation equation and the momentum equation
balancing the inertial force with the Cauchy stress 
with \eq{dafa-formula} and \eq{ultimate-}, we will thus deal with
the system for $\varrho$, $\vv$, $\Fe$, $\Lp$, and $\bm\xi$:
\begin{subequations}\label{Euler-hypoplast}
\begin{align}
\label{Euler0-hypoplast}&\DT\varrho=-\varrho\,{\rm div}\,\vv\,,\\
\nonumber
     &\varrho\DT\vv={\rm div}(\TT{+}\DD)
     +\varrho\GRAVITY\,\ \ \ \text{ with }\ \ 
\TT=\Se\Fee^{\!\!\top}\!\!+\varphi(\bm\xi,\Fe)\bbI\,,\ \ \ \Se=\varphi_{F}'(\bm\xi,\Fe)\,,
   \\[-.0em]
    &\hspace*{7.5em}\text{and }\ 
    \DD=\bbD\strain(\vv)-{\rm div}\,\mathfrak{H}\ \ \text{ with }\ 
\ \mathfrak{H}=\HYPER|\nabla^2\vv|^{p-2}\nabla^2\vv\,,
\label{Euler1-hypoplast}
\\[-.3em]\label{Euler2-hypoplast}
&\DT{\Fe}=
(\Nabla\vv)\Fe-\Fe\Lp\,,
\\\label{Euler3-hypoplast}
&\zetap(\bm\xi,\cdot)'(\Lp)={\rm dev}\MM\ \
\text{ with the Mandel stress }\
\MM=\Fee^{\!\!\top}\Se\,,\ \text{ and}
\\\label{Euler4-hypoplast}
&\DT{\bm\xi}=0\,,
\end{align}\end{subequations}
where $\GRAVITY$ denotes a prescribed acceleration in the actual frame
- most typically the gravity acceleration.
Moreover, the conservative Cauchy stress $\TT$ in \eq{Euler1-hypoplast}
has been additively enhanced by the {\it Stokes}-type {\it viscosity} stress
tensor $\DD$. For the slightly peculiar concept of the multipolar nonsimple material,
i.e.\ the so-called {\it hyperstress} $\mathfrak{H}$ contributing the linear Stokes
viscosity with the viscosity tensor $\bbD$ in \eq{Euler1-hypoplast}, we refer to
\cite{BeBlNe92PBMV,FriGur06TBBC,NeNoSi89GSIC,NecRuz92GSIV,PoGiVi10HHCS,Roub22VELS}.
This higher-order viscosity is used here primarily for analytical reasons to
facilitate convergence of the time-discrete scheme below and existence of
weak solutions of the system \eq{Euler-hypoplast}. Rather for notational simplicity,
we consider the linear Stokes-type viscosity using the tensor $\bbD$ and
(only scalar coefficient) $\HYPER$, both independent of $X$ and $F$.

The momentum equation \eq{Euler1-hypoplast} is to be supplemented by the
boundary conditions on the boundary $\varGamma$ of the domain $\varOmega$.
Due to the 2nd-grade multipolar concept, it is  a  4th-order
equation in space, so two boundary conditions are needed. 
Moreover, we consider the fixed shape of
$\varOmega$ by zero normal velocity while otherwise the boundary is free,
so we distinguish the normal and the tangential component. Altogether, we
consider
\begin{align}
&\vv\Cdot\nn=0,\ \ \ [(\TT{+}\DD)\nn+\divS(\mathfrak{H}\nn)]_\text{\sc t}^{}
=\bm0\,,\ \ \text{ and }\ \
\mathfrak{H}\Colon(\nn{\otimes}\nn)={\bm0}\ \ \text{ on }\ \varGamma\,,
\label{Euler-small-BC-hyper}\end{align}
 where $\divS={\rm tr}(\nablaS)$ is the surface divergence with
${\rm tr}(\cdot)$ denoting the trace of a 
2${\times}$2-matrix and with $\nablaS v$ denoting the surface gradient of
$v$ defined as $\nablaS v=\nabla v-\frac{\partial v}{\partial\nn}\nn$. 

Let us comment that fixing the shape of $\varOmega$ by setting $\vv\Cdot\nn=0$
is a very common simplification used in both theoretical and computational
studies in the Eulerian time-evolving frame.

The energetics of the model \eq{Euler-hypoplast} can be seen by testing
the momentum equation \eq{Euler1-hypoplast} by $\vv$ while using
the continuity equation \eq{Euler0-hypoplast} multiplied by $|\vv|^2/2$
and the kinematic equations \eq{Euler2-hypoplast} and \eq{Euler4-hypoplast}
multiplied respectively by $\varphi_{F}'(\bm\xi,\Fe)$ and
$\varphi_X'(\bm\xi,\Fe)$, and eventually using \eq{Euler3-hypoplast}
multiplied by $\Lp$. In particular, the test of the Cauchy stress $\TT$,
i.e.\ the term ${\rm div}\TT\Cdot\vv$, is treated by the Green formula as
\begin{align}\nonumber
\!\!\int_\varGamma\!&\vv\Cdot\TT\nn\,\d S-\!\int_\varOmega\!{\rm div}\TT\Cdot\vv\,\d\xx
=\!\int_\varOmega\!\big(\varphi_F'(\bm\xi,\Fe)\Fe^\top\big)\Colon\Nabla\vv
+\varphi(\bm\xi,\Fe)\,{\rm div}\,\vv\,\d\xx
\\[-.4em]&\nonumber
=\!\int_\varOmega\!
\varphi_F'(\bm\xi,\Fe)\Colon\big((\Nabla\vv)\Fe\big)+\varphi(\bm\xi,\Fe)\,{\rm div}\,\vv\,\d\xx
\\&\nonumber
=\!\int_\varOmega\!\varphi_F'(\bm\xi,\Fe)\Colon\big(\DT\Fe\!{+}
\Fe\Lp\big)+\varphi(\bm\xi,\Fe)\,{\rm div}\,\vv\,\d\xx
\\&\nonumber
=\int_\varOmega\!\varphi_F'(\bm\xi,\Fe)\Colon\DT\Fe+\MM\Colon\Lp+\varphi(\bm\xi,\Fe)\,{\rm div}\,\vv+\!\!\!\!\!\lineunder{\varphi_X'(\bm\xi,\Fe)\Cdot\DT{\bm\xi}}{=0 due to \eq{Euler4-hypoplast}\hspace{-2em}}\!\!\!\d\xx
\\[-.8em]&\nonumber
=\int_\varOmega\!\bigg(\varphi_X'(\bm\xi,\Fe)\Colon\pdt{\bm\xi}
+\varphi_F'(\bm\xi,\Fe)\Colon\pdt{\Fe}+\MM\Colon\Lp
\\[-.8em]&\nonumber\hspace{6em}
+\!\!\!\!\lineunder{\varphi_X'(\bm\xi,\Fe)\Colon(\vv\Cdot\Nabla)\bm\xi\!
+\varphi_F'(\bm\xi,\Fe)\Colon(\vv\Cdot\Nabla)\Fe\!
+\varphi(\bm\xi,\Fe)\,{\rm div}\,\vv}
{$\hspace{5em}=\int_\varGamma\varphi(\bm\xi,\Fe)(\vv\Cdot\nn)\,\d S=0
\text{ due to \eq{Euler-small-BC-hyper}}\hspace{-2em}$}\!\!\!\bigg)\,\d\xx
\\[-.9em]&
=\frac{\d}{\d t}\!\int_\varOmega\!\varphi(\bm\xi,\Fe)\,\d\xx+\!\int_\varOmega\!
\MM\Colon\Lp\,\d\xx\,.
\label{Euler-hypoplast-test-momentum}\end{align}

The boundary conditions \eq{Euler-small-BC-hyper} are also used for
the dissipative term ${\rm div}\DD\Cdot\vv$ while the inertial term
$\varrho\DT\vv\Cdot\vv$ leads to the kinetic energy $\frac12\varrho|\vv|^2$,
see \cite{Roub25TDVE} for details in calculations in particular as far as the
hyperstress and the higher-order boundary conditions \eq{Euler-small-BC-hyper}
concerns. Thus we finally obtain the {\it energy-dissipation balance}
\begin{align}
\!\!\!\!\!\frac{\d}{\d t}\int_\varOmega\!\!\!\!
\linesunder{\frac\varrho2|\vv|^2}{kinetic}{energy}
\!\!\!\!\!+\!\!\!\!\!\linesunder{\varphi(\bm\xi,\Fe)}{stored}{energy}\!\!\!\!
\,\d\xx
+\!\int_\varOmega\hspace{-.7em}
\linesunder{\bbD\strain(\vv)\Colon\strain(\vv){+}\HYPER|\nabla^2\vv|^p}{dissipation rate due to}{the Stokes-type viscosity}\hspace{-.7em}
+\hspace{-1em}\threelinesunder{\Lp\Colon
\MM}{dissipation rate}{due to  Maxwell}{viscosity}\hspace{-1em}\d\xx
=\int_\varOmega\hspace{-1.4em}\linesunder{\varrho\GRAVITY\Cdot\vv}{\ power of}{gravity field\ }\hspace{-1.4em} \,\d\xx\,.\ 
\nonumber\\[-1.5em]\label{EUL-L-engr-balance}\end{align}

In contrast to \cite{Roub25TDVE}, we avoid a gradient of the inelastic
distortion rate, which  (beside  simplification of the model) 
still allows for a rigorous analysis under a suitable qualification of
this viscoplastic dissipation potential, cf.\ \eq{EUL-L-ass-zeta} below.

\bigskip

\section{The staggered time discretization}\label{sec-discretization}

As in \cite{Roub25TDVE},
we express the system in terms of the linear momentum $\pp=\varrho\vv$ and,
later, exploit the fact that the kinetic energy $(\pp,\varrho)\mapsto
\frac12|\pp|^2/\varrho$ is convex, in contrast to the equivalent form
$(\vv,\varrho)\mapsto\frac12\varrho|\vv|^2$ which is nonconvex.
Moreover, we merge \eq{Euler2-hypoplast} with \eq{Euler3-hypoplast}
by writing more explicitly $\Lp=[\zetap(\bm\xi,\cdot)']^{-1}({\rm dev}\MM)$
and using the construction of the convex conjugate and the calculus from
the convex analysis for $[\zetap(X,\cdot)']^{-1}=[\zetap(X,\cdot)^*]'$, where
\begin{align}
\zetap(X,\cdot)^*(M)=\mbox{$\sup_{L\in\Rdev}$}\big(L\Colon M-\zetap(X,L)\big)
\end{align}
denotes the convex conjugate of the convex function $\zetap(X,\cdot)$.
Actually, writing the flow rule \eq{Euler3-hypoplast} in terms of
the convex conjugate in \eq{Euler2-hypoplast-heter} may have an advantage
if several Maxwellian-type visco-plastic mechanisms are considered in
series, as often in engineering or geophysical models, because
the convex conjugate is then easy to calculate explicitly in contrast
to the original potential $\zetap$, cf.\ \cite{Roub25FNVR}.
The system \eq{Euler-hypoplast} then reduces and modifies to a  system 
for $(\varrho,\pp,\Fe,\bm\xi)$ and thus also for $\vv=\pp/\varrho$ as
\begin{subequations}\label{Euler-hypoplast-p-heter}
\begin{align}
\label{Euler0-hypoplast-heter}&\pdt\varrho=-{\rm div}\,\pp\ \ \ \text{ with }\ 
\pp=\varrho\vv\,,\\
\nonumber
     &\pdt\pp={\rm div}(\TT{+}\DD{-}\pp{\otimes}\vv)
     +\varrho\GRAVITY\,\ \ \ \text{ with }\ \ 
\TT=\Se\Fee^{\!\!\top}\!\!+\varphi(\bm\xi,\Fe)\bbI\,,\ \ \ \Se=\varphi_{F}'(\bm\xi,\Fe)\,,
   \\[-.3em]
    &\hspace*{7.2em}\text{and }\ 
    \DD=\bbD\strain(\vv)-{\rm div}\,\mathfrak{H}\ \ \text{ with }\ 
\ \mathfrak{H}=\HYPER|\nabla^2\vv|^{p-2}\nabla^2\vv\,,
\label{Euler1-hypoplast-heter}
\\\label{Euler2-hypoplast-heter}
&\pdt{\Fe}=
(\Nabla\vv)\Fe-(\vv\Cdot\nabla)\Fe
-\Fe\big[\zetap(\bm\xi,\cdot)^*\big]'\big({\rm dev}(\Fee^{\!\!\top}\Se)\big)\,,
\\\label{Euler3-hypoplast-heter}
&\pdt{\bm\xi}=-(\vv\Cdot\nabla)\bm\xi\,
\end{align}\end{subequations}

We will consider the initial-value problem for the system
\eq{Euler-hypoplast-p-heter} with the boundary conditions
\eq{Euler-small-BC-hyper}
by prescribing 
\begin{align}\label{IC-large}
\varrho|_{t=0}^{}=\varrho_0^{}\,,\ \ \ \ 
\pp|_{t=0}^{}=\varrho_0^{}\vv_{0}^{}\,,\ \ \ \ {\Fe}|_{t=0}^{}=\Fezero\,,
\ \text{ and }\ {\bm\xi}|_{t=0}^{}=\text{\,identity on }\,\varOmega\,,
\end{align}
and analyze the model \eq{Euler-hypoplast} in terms of
its weak solutions, introducing:

\begin{definition}[Weak formulation of \eq{Euler-hypoplast-p-heter}]\label{def-ED-Ch5}
We call $(\varrho,\vv,\Fe,\bm\xi)$ with $\varrho\in W^{1,1}(I{\times}\varOmega)$,
$\vv\in L^p(I;W^{2,p}(\varOmega;\R^3))$, $\Fe\in W^{1,1}(I{\times}\varOmega;\R^{3\times3})$,
and $\bm\xi\in W^{1,1}(I{\times}\varOmega;\R^3)$ a weak solution to the system
\eq{Euler-hypoplast} with the boundary conditions \eq{Euler-small-BC-hyper}
and the initial conditions \eq{IC-large} if 
\begin{align}
&\nonumber
\int_0^T\!\!\!\int_\varOmega\bigg(\Big(\varphi_F'(\bm\xi,\Fe)\Fe^\top\!+
\bbD\strain(\vv)-\varrho\vv{\otimes}\vv\Big){:}\strain(\widetilde\vv)
+\varphi(\bm\xi,\Fe)\,{\rm div}\,\widetilde\vv
-\varrho\vv{\cdot}\pdt{\widetilde\vv}
\\[-.6em]&\hspace{4em}
+\HYPER|\nabla^2\vv|^{p-2}\nabla^2\vv\Vdots\nabla^2\widetilde\vv\bigg)\,\d\xx\d t
=\!\int_0^T\!\!\!\int_\varOmega\varrho\GRAVITY{\cdot}\widetilde\vv\,\d\xx\d t
+\!\int_\varOmega\!\varrho_0^{}\vv_0^{}{\cdot}\widetilde\vv(0)\,\d\xx
\label{def-ED-Ch5-momentum}\end{align}
holds for any $\widetilde\vv$ smooth,  say $\widetilde\vv\in
L^\infty(I;W^{2,\infty}(\varOmega;\R^3))\,\cap\,W^{1,\infty}(I;L^\infty(\varOmega;\R^3))$ 
with $\widetilde\vv{\cdot}\nn={\bm0}$ on
$I{\times}\varGamma$ and $\widetilde\vv(T)=0$, and if the equations
(\ref{Euler-hypoplast-p-heter}a,c,d)
hold a.e.\ on $I{\times}\varOmega$ together with the respective initial conditions
for $\varrho$, $\Fe$, and ${\bm\xi}$ in \eq{IC-large}.
\end{definition}

\def\GL{{\rm GL}}

Let us summarize the assumptions used in what follows:
\begin{subequations}\label{EUL-L-ass}
\begin{align}\label{EUL-L-ass-phi}
&\!\varphi\in C^1(\barOmega{\times}\GL_3^+), \ \
\inf_{X\in\varOmega,\ F\in \GL_3^+}\varphi(X,F)>-\infty\,,\ \ \
\forall X{\in}\varOmega:\lim_{{\rm det}F\to0}\varphi(X,F)=+\infty\,,
\\[-.1em]&\nonumber
\!\zetap\in C(\varOmega{\times}\Rdev),\ \ \ \ 
\forall X{\in}\varOmega:\ \ \
\zetap(X,\cdot)\ \text{is convex with a minimum at }\{0\},\
\\&\nonumber
\qquad\qquad\forall f\in W^{1,\infty}(\varOmega{\times}\R^{3\times3};\R^{3\times3})\ \ \forall B\subset\GL_3^+\text{ compact}:
\\[-.1em]&\nonumber
\qquad\qquad\qquad\qquad\qquad\
(X,F)\mapsto F\big[\zetap(X,\cdot)^*\big]'\big({\rm dev}(F^\top\!f(X,F))\big)
\\[-.1em]&\nonumber
\qquad\qquad\qquad\qquad\qquad\
\text{ is Lipschitz continuous as a mapping $\ \varOmega{\times}B\!\to\R^{3\times3}$,
and}
\\&\label{EUL-L-ass-zeta}
\qquad\qquad 0<\!\!\inf_{X\in\varOmega,\ L\in \Rdev\setminus\{0\}}\!\!\!\!
\frac{\zetap(X,L)}{|L|^2}\ \ \text{ and }\ 
\sup_{X\in\varOmega,\ L\in \Rdev}\!\!\frac{\zetap(X,L)}{1{+}|L|^2}
<+\infty\,,
\\[-.3em]&\!\label{Euler-small-ass-D}
\bbD:\Rsym\to\Rsym\ \text{ linear symmetric positive semidefinite}\,,
\ \ \ \HYPER>0\,,
\\[-.1em]&\!
\vv_0\,{\in}\,L^2(\varOmega;\R^3),\ \ \varrho_0\,{\in}\,W^{1,r}(\varOmega),\ \ 
\min_{\barOmega}\varrho_0>0,\ \
\Fezero\,{\in}\,W^{1,r}(\varOmega;\R^{3\times3}),
\ \min_{\barOmega}\det\Fezero>0,\!\label{EUL-L-ass-IC}
\\[-.5em]&\!\label{Euler-small-ass-g}
\GRAVITY\in L^\infty(I{\times}\varOmega;\R^3)\,. 
\end{align}\end{subequations}
 The applicability of (\ref{EUL-L-ass}a,b) is discussed
in Remark~\ref{rem-ass} below.

To devise the stable decoupled time discretization, let us first introduce
a truncated stored energy $\varphi_\lambda$ which vanishes for too large $F$
or too small $\det F$. Specifically,
\begin{align}\nonumber
\varphi_\lambda(X,F)=\begin{cases}\qquad\ \varphi(X,F)&
\text{for $|F|\le \lambda$ and $\det F\ge\Frac1\lambda$},
\\[.3em]\qquad\qquad0&\text{for $|F|\ge 2\lambda$ or $\det F\le \Frac1{2\lambda}$},
\\\Big(3\lambda^2\big(2\det F{-}\Frac 1\lambda\big)^2-2\lambda^3\big(2\det F{-}\Frac 1\lambda\big)^3\Big)\times\hspace*{-2.5em}&
\\\qquad
\times\Big(3\big(\Frac{|F|}\lambda{-}1\big)^2-2\big(\Frac{|F|}\lambda{-}1\big)^3\Big)
\varphi(X,F)&\text{otherwise}\,,
\end{cases}\!\!\!\!\!\!
\\[-2.1em]\label{cutoff-stress}
\end{align}
where $|\cdot|$ is the Frobenius norm of a matrix, i.e.\ $|F|=(\sum_{i,j=1}^3
F_{ij}^2)^{1/2}$, which makes $\varphi_\lambda(X,\cdot)$ frame indifferent
if $\varphi(X,\cdot)$ is so; cf.\
\cite{Roub24TVSE,RouSte23VESS,RouTom23IFST} for this truncation.

Considering the equidistant time partition of the time interval
$I=[0,T]$ with the  time  step $\tau>0$ with $T/\tau$ integer, the
proposed staggered scheme reads as 
\begin{subequations}\label{EUL-L-disc}
\begin{align}\label{EUL-L-disc0}
&\!\!\frac{\varrho_\tau^k{-}\varrho_\tau^{k-1}\!\!}\tau\,=-
{\rm div}\,\pp_\tau^k\ \ 
\text{ with }\,\pp_\tau^k=\varrho_\tau^k\vv_\tau^k\,,\!
\\[-.2em]&\nonumber
\!\!\frac{\pp_\tau^k{-}\pp_\tau^{k-1}\!\!}\tau\,=
{\rm div}\Big(
\mathscr{T}_\lambda(\bm\xi_{\tau}^{k-1},\Fekk)
{+}\DD_\tau^k{-}\pp_\tau^k{\otimes}\vv_\tau^k\Big)
+\varrho_\tau^k\GRAVITY_{\tau}^k\,,
\\[-.2em]&
\hspace*{5em}\text{ where }\ \DD_\tau^k=\bbD\strain(\vv_\tau^k)
-{\rm div}\,\mathfrak{H}_\tau^k\ \ \text{ with }\ \
\mathfrak{H}_\tau^k=\HYPER\big|\nabla^2\vv_\tau^k\big|^{p-2}\nabla^2\vv_\tau^k\,,\!
\label{EUL-L-disc1}
\\[-.1em]
&\!\!\frac{\Fek{-}\Fekk\!\!}\tau\,
=(\nabla\vv_\tau^k)\Fek-
\Fek\mathscr{L}_\lambda^{}(\bm\xi_{\tau}^{k},\Fek)
-(\vv_\tau^k\Cdot\nabla)\Fek\,,\ \ \text{ and}
\label{EUL-L-disc2}
\\[-.1em]
&\!\!\frac{\bm\xi_{\tau}^k{-}\bm\xi_{\tau}^{k-1}\!\!}\tau\,
=-(\vv_\tau^k\Cdot\nabla)\bm\xi_{\tau}^k\,
\label{EUL-L-disc3}
\end{align}\end{subequations}
with the Cauchy stress and the inelastic distortion rate truncated
\begin{subequations}\label{truncated}\begin{align}\label{stress-truncated+}
&\mathscr{T}_\lambda^{}(X,F):=[\varphi_\lambda]_F'(X,F)F^\top\!{+}\varphi_\lambda^{}(X,F)\,\bbI
\ \ \ \text{ and }\  
\\&\label{Lp-truncated-heter}
\mathscr{L}_\lambda^{}(X,F):=\big[\zetap(X,\cdot)^*\big]'\big({\rm dev}(F^\top\![\varphi_\lambda^{}(X,\cdot)]'(F))\big)
\end{align}\end{subequations}
with $\varphi_\lambda$ from \eq{cutoff-stress}. Moreover,
$\GRAVITY_{\tau}^k(\xx)=\frac1\tau\int_{(k-1)\tau}^{k\tau}\GRAVITY(t,\xx)\,\d t$.
The index $k$ ranges $1,2,...,T/\tau$, and the solution \eq{EUL-L-disc}
is to approximate solution to \eq{Euler-hypoplast-p-heter} at time $t=k\tau$.
The corresponding boundary conditions for \eq{EUL-L-disc1} are now
\begin{align}\nonumber
&\vv_\etau^k{\cdot}\nn=0\,,\ \ \ \ \ \nabla^2\vvk{:}(\nn{\otimes}\nn)=0\,,
\ \text{ and }\
\\&\Big[\Big(\mathscr{T}_\lambda(\bm\xi_{\etau}^{k-1},\Fekk)
{+}\bbD\strain(\vvk){-}{\rm div}\mathfrak{H}_\etau^k\Big)\nn
-\divS\big(\mathfrak{H}_\etau^k\Cdot\nn\big)\Big]_\text{\sc t}^{}\!\!=\bm0\,.
\label{BC-disc}\end{align}

Note that the system \eq{EUL-L-disc} is indeed decoupled (staggered) in
the sense that (\ref{EUL-L-disc}a,b) gives
$(\varrho_\etau^k,\pp_\etau^k)$ and thus also $\vv_\etau^k$. This $\vv_\etau^k$ is
then used for \eq{EUL-L-disc3} to obtain $\bm\xi_\etau^k$ and eventually
for \eq{EUL-L-disc2} to obtain $\Fek$.

Notably,  unlike in \cite{Roub25TDVE}, the approximation scheme is cast
without any higher gradients into the kinematic equations
and, perhaps more importantly, no regularization of the model is needed and thus
a direct (not successive) limit passage within the time step $\tau\to0$
can be done.
Moreover, unlike in \cite{Roub25TDVE}, 
the general (not necessarily polyconvex) stored energy  and avoid a separate
kinematics of the minors (in particular of $\det\Fe$ to be discretized,
and 
we will  use the
discrete Gronwall inequality and a sparsity estimate in a careful way.

We will  use  the notation that exploits the interpolants defined as
follows: using the values $(\bm\xi_\etau^k)_{k=0}^{T/\tau}$,
we define the piecewise constant and the piecewise  affine interpolants
respectively as
\begin{align}\nonumber
&\overlinexitau(t)\!:=\bm\xi_\etau^k\,,\ \ \ 
\underlinexitau(t)\!:=\bm\xi_\etau^{k-1}\,,
\ \ \text{ and }\ \ 
\\&\bm\xi_\etau(t)\!:=\Big(\frac t\tau{-}k{+}1\Big)\bm\xi_\etau^k
\!+\Big(k{-}\frac t\tau\Big)\bm\xi_\etau^{k-1}\
\ \ \text{ for }\ (k{-}1)\tau<t\le k\tau
\label{def-of-interpolants}
\end{align}
for $k=0,1,...,T/\tau$. Analogously, we define also
$\overlineEetau$, $\underlineFetau$, $\Fe_{,\tau}$, $\overlinerhoetau$,
$\varrho_\etau^{}$, $\overlineppetau$, $\pp_\etau^{}$, etc.
In terms of these interpolants, the recursive system \eq{EUL-L-disc}
can be written ``more compactly'' as
\begin{subequations}\label{EUL-L-visco-ED+discr}
\begin{align}\label{EUL-L-visco-ED+0discr}
&\!\!\pdt{\varrho_\etau}=
-{\rm div}\,\overlineppetau\ \
\text{ with $\ {\overlineppetau}=\overlinerhoetau\overlinevvtau$}\,,
\\&\nonumber
\!\!\pdt{\pp_\etau}=
{\rm div}\Big(\mathscr{T}_\lambda(\underlinexitau,\underlineFetau)
+\overlineDetau-
\overlineppetau{\otimes}\overlinevvtau\Big)
+\overlinerhoetau\overline\GRAVITY_{\tau}\,,
\\[-.2em]&\hspace*{3em}\text{ where }
\ \overlineDetau=\bbD\strain(\overlinevvtau)
-{\rm div}\overline{\mathfrak{H}}_\etau\ 
\text{ with }\ \overline{\mathfrak{H}}_\etau=
\HYPER\big|\nabla^2\overlinevvtau\big|^{p-2}\nabla^2\overlinevvtau\,,
\label{EUL-L-visco-ED+1discr}
\\[-.3em]&
\!\!\pdt{\Fe_{,\tau}}
=(\nabla\overlinevvtau)\overlineEetau
-\overlineEetau\mathscr{L}_\lambda(\overlinexitau,\overlineEetau)
-(\overlinevvtau\Cdot\nabla)\overlineEetau\,,\ \ \text{ and}
\label{EUL-L-visco-ED+2discr}
\\[-.0em]&
\!\!\pdt{\bm\xi_\etau}
=-(\overlinevvtau\Cdot\nabla)\overlinexitau
\label{EUL-L-visco-ED+3discr}
\end{align}\end{subequations}
with the corresponding boundary conditions and with the initial conditions
for $(\varrho_\etau,\pp_\etau,\Fe_{,\tau},\bm\xi_\etau)$.

\begin{proposition}[Stability and convergence of \eq{EUL-L-visco-ED+discr}
and existence of solutions to \eq{Euler-hypoplast-p-heter}]
\label{prop-ED-Ch5-existence}
$_{}$\linebreak
Let \eq{EUL-L-ass}  hold  with $ p> r>3$. Then:\\
\ITEM{(i)}{For any $\lambda$ and for all sufficiently small time steps
$\tau>0$, the staggered scheme \eq{EUL-L-disc} has a solution
$(\varrho_\etau^k,\vv_\etau^k,\Fek,\bm\xi_\etau^k)\in W^{1,r}(\varOmega)
\times W^{2,p}(\varOmega;\R^{3})\times W^{1,r}(\varOmega;\R^{3\times3})\times
W^{1, r_1}(\varOmega;\R^3)$ and $\varrho_\etau^k>0$ on  $\barOmega$.
}
\ITEM{(ii)}{For any $\lambda$ fixed and for $\tau\to0$, there is a subsequence
of $(\overlinerhoetau,\overlinevvtau,\overlineEetau,\overlinexitau)$ converging
in the topologies specified in
\eq{Euler-small-converge} below to some $(\varrho,\vv,\Fe,\bm\xi)$. Moreover,
there is a sufficiently large $\lambda$ for which any such quadruple is a weak
solution in the sense of Definition~\ref{def-ED-Ch5}.}
\ITEM{(iii)}{In particular, the initial-boundary-value
problem for the system \eq{Euler-hypoplast-p-heter} has at least one weak
solution $(\varrho,\vv,\Fe,\bm\xi)$ such that
also $\varrho\in L^\infty(I;W^{1,r}(\varOmega))\,\cap\,C(I{\times}\barOmega)$ with
$\min_{I{\times}\barOmega}\varrho>0$, $\vv\in L^\infty (I;L^2(\varOmega;\R^3))$,
$\Fe\in L^\infty(I;W^{1,r}(\varOmega;\R^{3\times3}))
\,\cap\,C(I{\times}\barOmega;\R^{3\times3})$ with $\min_{I\times\barOmega}\det\Fe>0$
with $r$ from \eq{EUL-L-ass-IC}, and
$\bm\xi\in L^\infty(I;W^{1, r_1}(\varOmega;\R^3))\,\cap\,C(I{\times}\barOmega;\R^3)$
with $ r_1<+\infty$ arbitrary. Moreover, the energy-dissipation balance
\eq{EUL-L-engr-balance} integrated over the time interval $[0,t]$ holds
for any $t\in I$.}
\end{proposition}

\begin{proof}[Sketch of the proof]
For clarity, we divide the proof into six steps. Some arguments are 
only sketched while referring to \cite{Roub25TDVE} which employed a similar 
(although regularized and fully implicit) time discretization of the
subsystem (\ref{EUL-L-disc}a,b) with a similar analysis
as far as the mass density and the velocity concerns.

\medskip{\it Step 1: Formal a-priori estimates and the choice of $\lambda$}.
We use (formally at this point) the energy-dissipation balance
\eq{EUL-L-engr-balance} together with $\varrho\ge0$.
Treating the gravitational loading as in
\eq{Euler-small-est-Gronwall} below, we obtain the a-priori bounds
$\|\strain(\vv)\|_{L^2(I{\times}\varOmega;\R^{3\times3})}^{}\le C$ and 
$\|\nabla^2\vv\|_{L^p(I{\times}\varOmega;\R^{3\times3\times3})}^{}\le C$.
This a-priori quality of the velocity field together with the qualification of 
the initial condition $\Fezero$ in \eq{EUL-L-ass-IC} allows for 
using the transport-and-evolution equation \eq{Euler2-hypoplast-heter}
to obtain the a-priori bound
$\|\Fe\|_{L^\infty(I;W^{1,r}(\varOmega;\R^{3\times3}))}\le C_r^{(1)}$.
For details see \cite[Lemma\,5.1]{Roub24TVSE}; actually, rather for analytical
arguments we need $\vv$ also in $L^1(I{\times}\varOmega;\R^3)$ although
not a-priori bounded, which is also seen from \eq{EUL-L-engr-balance}
provided $\min_{I{\times}\bar{\varOmega}}\varrho>0$ which
can be seen from \eq{Euler-hypoplast-p-heter} written for the
sparsity $\sigma:=1/\varrho$, cf.\ \eq{bound-for-sparsity} below.

Also, from \eq{Euler2-hypoplast-heter} together with isochoricity of the
inelastic distortion $\Fp$, we have also the equation for $1/\det\Fe$ as
$(1/\det\Fe)\!\DT{^{}\!}
=({\rm div}\,\vv)/\det\Fe$, cf.\
\eq{DT-det-extended+} below. From \eq{EUL-L-ass-IC}, also
$$
\nabla\frac1{\det\Fezero}=-\frac{\nabla(\det\Fezero)}{(\det\Fezero)^2}
=\Frac{{\rm Cof}\Fezero}{(\det\Fezero)^2}\nabla\Fezero\in L^r(\varOmega;\R^3)\,,
$$
so that we can infer $\|1/\det\Fe\|_{L^\infty(I;W^{1,r}(\varOmega))}\le C_r^{(2)}$.

This suggests to put $\lambda>N_r\max(C_r^{(1)},C_r^{(2)})$ with $N_r$ denoting
the norm of the embedding $W^{1,r}(\varOmega)\subset L^\infty(\varOmega)$
holding for $r>3$ from the condition for $\Fezero$ in \eq{EUL-L-ass-IC}.
 We will exploit this choice of $\lambda$ sufficiently large
in Steps~5-6, while Steps~2-4 work for arbitrary $\lambda\ge1$.
It should be emphasized that {\it no limit passage} for $\lambda\to\infty$
will be needed if $\lambda$ will be considered large as specified here. 

\medskip{\it Step 2: Basic stability of the scheme
\eq{EUL-L-visco-ED+discr} and first a-priori estimates}.
Like used for \eq{EUL-L-engr-balance}, 
testing \eq{EUL-L-disc1} by $\vv_\etau^k$ and
using \eq{EUL-L-disc0} tested by $|\vv_\etau^k|^2/2$,
after summation over time steps from 1 up to $k$, we obtain
\begin{align}\nonumber
\!\!\!\!\!&\int_\varOmega\!\frac{|\pp_\etau^k|^2\!\!\!}{2\varrho_\etau^k\!\!\!}\,\d\xx
+\tau\!\sum_{m=1}^k\int_\varOmega\!
\bbD\strain(\vv_\etau^m)\Colon\strain(\vv_\etau^m)
+\HYPER|\nabla^2\vv_\etau^m|^p\,\d\xx
\\[-.5em]&\hspace{1.3em}
\le\int_\varOmega\!\frac{|\pp_0|^2\!\!}{2\varrho_0\!\!}\,\d\xx
+\tau\!\sum_{m=1}^k\int_\varOmega\!\varrho_\etau^m\GRAVITY_{\tau}^m\Cdot\vv_\etau^m
-\mathscr{T}_\lambda^{}(\bm\xi_{\eetau}^{m-1},\Femm)\Colon\strain(\vv_\etau^m)
\,\d\xx\!\!
\label{EUL-L-basic-engr-balance-disc}\end{align}
Since $\mathscr{T}_\lambda^{}(\bm\xi_{\eetau}^{m-1},\FF_{\eetau}^{m-1})$ is a-priori
bounded in $L^\infty(\varOmega;\R^{3\times 3})$,
the last term in \eq{EUL-L-basic-engr-balance-disc} can easily be
treated by the assumed positive definiteness \eq{Euler-small-ass-D}
of $\bbD$ in the left-hand side, while the penultimate term
$\varrho_\etau^m\GRAVITY_{\tau}^m\Cdot\vv_\etau^m$ can be estimated as
\begin{align}
\nonumber
\!\!\int_\varOmega&\varrho_\etau^m\GRAVITY_{\tau}^m\Cdot\vv_\etau^m\,\d\xx
=\int_\varOmega\pp_\etau^m\Cdot\GRAVITY_{\tau}^m\,\d\xx\le
\Bigg\|\frac{\pp_\etau^m}{\sqrt{\varrho_\etau^m}}\Bigg\|_{L^2(\varOmega;\R^3)}^{}\!\!
\big\|\!\sqrt{\varrho_\etau^m}\GRAVITY_{\tau}^m\big\|_{L^2(\varOmega;\R^3)}^{}
\\[-.6em]&
\le\Bigg\|\frac{\pp_\etau^m}{\sqrt{\varrho_\etau^m}}\Bigg\|_{L^2(\varOmega;\R^3)}^{}\!\!
\|\varrho_\etau^m\|_{L^1(\varOmega)}^{}\|\GRAVITY_{\tau}^m\|_{L^\infty(\varOmega;\R^3)}^{}
\le \sqrt R\|\GRAVITY_{\tau}^m\|_{L^\infty(\varOmega;\R^3)}^{}
\Bigg(\!1+
\!\Frac12\!\int_\varOmega\!\frac{|\pp_\etau^{ m}|^2\!\!\!}{2\varrho_\etau^{ m}\!\!\!}
\,\d\xx\Bigg)\,.
\label{Euler-small-est-Gronwall}\end{align}
with $R=\int_\varOmega\varrho_0\,\d\xx$; here we used the mass conservation 
guaranteed by \eq{EUL-L-disc0} together with the first boundary
condition in \eq{BC-disc}. Then, for
$\tau<1/(\sqrt R\|\GRAVITY\|_{L^\infty(I\times\varOmega;\R^3)}^{})$,
we can use the discrete Gronwall inequality \eq{9.2++} to obtain the a-priori
estimates
\begin{subequations}\label{Euler-small-est}\begin{align}
\label{Euler-small-est1}
&\bigg\|\frac{\overlineppetau}{\sqrt{\overlinerhoetau}}
\bigg\|_{L^\infty(I;L^2(\varOmega;\R^3))}^{}\!\!\le C\,,\ \ \ 
\|\strain(\overlinevvtau)\|_{L^2(I{\times}\varOmega;\R^{3\times3})}^{}\le C\,,\ \ \
\|\nabla^2\overlinevvtau\|_{L^p(I{\times}\varOmega;\R^{3\times3\times3})}^{}\le C
\intertext{and then, from the first estimate in \eq{Euler-small-est1},
we have also}
&\|\overlineppetau\|_{L^\infty(I;L^1(\varOmega;\R^3))}^{}\le\!\!\!\!
\lineunder{\big\|\!\sqrt{\overlinerhoetau}\big\|_{L^\infty(I;L^2(\varOmega))}^{}}
{$\ \ =\sqrt R$}\!\!\!
\bigg\|\frac{\overlineppetau}{\sqrt{\overlinerhoetau}}\bigg\|_{L^\infty(I;L^2(\varOmega;\R^3))}^{}\!\!
\le C
\label{Euler-small-est3+++}\end{align}\end{subequations}
with $C$ here and below denoting a generic constant.

Then \eq{EUL-L-disc0} can be tested by 
$|\varrho_{\etau}^k|^{s-2}\varrho_{\etau}^k$ with some $s>1$. Using the Green
formula with the boundary condition $\nn\Cdot\vv_{\etau}^k=0$, the convective
term can be handled as
\begin{align}\nonumber
\!\!\int_\varOmega&{\rm div}(\varrho_{\etau}^k\vv_{\etau}^k)
(|\varrho_{\etau}^k|^{s-2}\varrho_{\etau}^k)\,\d\xx
=\!\!\int_\varOmega({\rm div}\,\vv_{\etau}^k)|\varrho_{\etau}^k|^s
+(\vv_{\etau}^k\Cdot\nabla\varrho_{\etau}^k)
(|\varrho_{\etau}^k|^{s-2}\varrho_{\etau}^k)\,\d\xx
\\[-.3em]&
=\!\int_\varOmega\!({\rm div}\,\vv_{\etau}^k)|\varrho_{\etau}^k|^s
-\varrho_{\etau}^k\,{\rm div}\big(|\varrho_{\etau}^k|^{s-2}\varrho_{\etau}^k\vv_{\etau}^k\big)\,\d\xx
=\Big(1{-}\Frac1s\Big)\int_\varOmega({\rm div}\,\vv_{\etau}^k)|\varrho_{\etau}^k|^s\d\xx\,,
\label{transport-rho}\end{align}
so that this test gives the inequality
\begin{align}
  &\int_\varOmega\!\frac{|\varrho_{\etau}^k|^s\!-|\varrho_\tau^{k-1}|^s\!}\tau
  \,\d\xx\,\le\!\Frac1{s'}\int_\varOmega\!({\rm div}\,\vv_\tau^k)|\varrho_{\etau}^k|^s
\,\d\xx
\le\Frac1{s'}\|{\rm div}\,\vv_\tau^k\|_{L^\infty(\varOmega)}^{}\int_\varOmega\!|\varrho_{\etau}^k|^s\,\d\xx\,.
\label{transport-rho+}\end{align}
Due to \eq{Euler-small-est1} with $p>3$, we have
$\tau\sum_{k=1}^{T/\tau}\|{\rm div}\,\vv_\tau^k\|_{L^\infty(\varOmega)}^{ 2}\le C$,
from which we can see
$\max_{1\le k\le T/\tau}\|{\rm div}\,\vv_\tau^k\|_{L^\infty(\varOmega)}^{}\le(C/\tau)^{1/{ 2}}$.
This allows us to use the discrete Gronwall inequality \eq{9.2++} which,
for a{}ny  sufficiently small $\tau$, say for $\tau\le \frac12C$, gives
the estimate
\begin{align}\label{EUL-L-est-rho-s}
\|\overlinerhoetau\|_{L^\infty(I;L^s(\varOmega;\R^3))}^{}\le C\,.
\end{align}
Furthermore, we test \eq{EUL-L-disc0} by 
${\rm div}(|\nabla\varrho_{\etau}^k|^{r-2}\nabla\varrho_{\etau}^k)$ with some
$r>1$. Using the Green formula with the boundary condition
$\nn\Cdot\vv_{\etau}^k=0$, the convective term can be handled as
\begin{align}\nonumber
&\!\!\int_\varOmega\!\nabla\big(\vv_{\etau}^k\Cdot\nabla\varrho_{\etau}^k
  \big)\Cdot\big(|\nabla\varrho_{\etau}^k|^{r-2}\nabla\varrho_{\etau}^k\big)\,\d\xx
\\[-.5em]&\hspace{.1em}\nonumber
=\!\int_\varOmega|\nabla\varrho_{\etau}^k|^{r-2}(\nabla\varrho_{\etau}^k{\otimes}\nabla\varrho_{\etau}^k)\Colon\strain(\vv_{\etau}^k)
+(\vv_{\etau}^k\Cdot\nabla)\nabla\varrho_{\etau}^k
\Cdot\big(|\nabla\varrho_{\etau}^k|^{r-2}\nabla\varrho_{\etau}^k\big)\,\d\xx
\\[-.1em]&\hspace{.1em}\nonumber
=\!\int_\varGamma|\nabla\varrho_{\etau}^k|^r\,\vv_{\etau}^k\Cdot\nn\,\d S
+\!\int_\varOmega\bigg(|\nabla\varrho_{\etau}^k|^{r-2}(\nabla\varrho_{\etau}^k{\otimes}\nabla\varrho_{\etau}^k)\Colon\strain(\vv_{\etau}^k)
\\[-.6em]&\hspace{5.5em}\nonumber
-({\rm div}\,\vv_{\etau}^k)|\nabla\varrho_{\etau}^k|^r-(r{-}1)|\nabla\varrho_{\etau}^k|^{r-2}|\nabla\varrho_{\etau}^k|^{r-2}\nabla\varrho_{\etau}^k\Cdot
(\vv_{\etau}^k\Cdot\nabla)\nabla\varrho_{\etau}^k\bigg)\,\d\xx
\\[-.1em]&\hspace{.1em}
=\!\int_\varGamma\Frac1r|\nabla\varrho_{\etau}^k|^r\!\!\!\!
\lineunder{\vv_{\etau}^k\Cdot\nn}{$=0$}\!\!\!\!\d S
+\!\int_\varOmega|\nabla\varrho_{\etau}^k|^{r-2}(\nabla\varrho_{\etau}^k{\otimes}\nabla\varrho_{\etau}^k)\Colon\strain(\vv_{\etau}^k)-\Frac1r({\rm div}\,\vv_{\etau}^k)|\nabla\varrho_{\etau}^k|^r\,\d\xx\,.\!\!
\label{calulus-nonlin-rho-r}\end{align}
Thus, counting also with
$$
\int_\varOmega\nabla\big(({\rm div}\,\vv_{\etau}^k)\varrho_{\etau}^k\big)\Cdot
\big(|\nabla\varrho_{\etau}^k|^{r-2}\nabla\varrho_{\etau}^k\big)\,\d\xx=\!
\int_\varOmega({\rm div}\,\vv_{\etau}^k)|\nabla\varrho_{\etau}^k|^r
+\varrho_{\etau}^k|\nabla\varrho_{\etau}^k|^{r-2}\nabla\varrho_{\etau}^k
\Cdot\nabla({\rm div}\,\vv_{\etau}^k)\,\d\xx\,,
$$
we obtain the inequality
\begin{align}\nonumber
&\Frac1r
\int_\varOmega\frac{|\nabla\varrho_{\etau}^k|^r{-}|\nabla\varrho_{\etau}^{k-1}|^r
\!\!}\tau\,\d\xx
\le\!\int_\varOmega\bigg(\Big(\Frac1r{-}1\Big)
({\rm div}\,\vv_{\etau}^k)|\nabla\varrho_{\etau}^k|^r\!
-|\nabla\varrho_{\etau}^k|^{r-2}(\nabla\varrho_{\etau}^k{\otimes}
\nabla\varrho_{\etau}^k)\Colon\strain(\vv_{\etau}^k)
\\[-.3em]&\nonumber\qquad\qquad
-\varrho_{\etau}^k|\nabla\varrho_{\etau}^k|^{r-2}\nabla\varrho_{\etau}^k
\Cdot\nabla({\rm div}\,\vv_{\etau}^k)\bigg)\,\d\xx
\le 2
\big\|\strain(\vv_{\etau}^k)\big\|_{L^\infty(\varOmega;\R^{3\times3})}^{}
\int_\varOmega|\nabla\varrho_{\etau}^k|^r\,\d\xx
\\[-.3em]&\hspace{12em}+\|\varrho_{\etau}^k\|_{L^{ pr/(p-r)}(\varOmega)}^{}
\|\nabla({\rm div}\,\vv_{\etau}^k)\|_{L^p(\varOmega;\R^3)}^{}
\|\nabla\varrho_{\etau}^k\|_{L^r(\varOmega;\R^3)}^{r-1}
\,.\!\!
\label{calulus-nonlin-rho-r+}\end{align}
Using it with \eq{EUL-L-est-rho-s} for $s$  sufficiently
large, namely $s\ge pr/(p{-}r)$ exploiting here the assumption $r<p$ 
and using  
again the discrete Gronwall inequality \eq{9.2++} relying on
\eq{Euler-small-est1}, for  any  sufficiently small time steps $\tau$, 
we obtain the estimate
\begin{align}\label{est-of-rho-disc}
\|\overlinerhoetau\|_{L^\infty(I;W^{1,r}(\varOmega))}^{}\le C\,;
\end{align}
specifically, $\tau$'s should be small as it was for also for
\eq{EUL-L-est-rho-s} and in particular also $\nu>0$ assumed in
\eq{Euler-small-ass-D} is used. 

 Moreover,  \eq{EUL-L-disc3} is a vectorial but simpler
variant of \eq{EUL-L-disc0}. Thus, when tested \eq{EUL-L-disc3} by 
$|\bm\xi_{\etau}^k|^{s-2}\bm\xi_{\etau}^k$ with some $s>1$,
we obtain the inequality
\begin{align}
  &\int_\varOmega\!\frac{|\bm\xi_\tau^k|^s\!-|\bm\xi_\tau^{k-1}|^s\!}\tau
  \,\d\xx\,\le\!\int_\varOmega\!({\rm div}\,\vv_\tau^k)|\bm\xi_\tau^k|^s\,\d\xx
\le\|{\rm div}\,\vv_\tau^k\|_{L^\infty(\varOmega)}^{}\int_\varOmega\!|\bm\xi_\tau^k|^s\,\d\xx\,
\label{transport+}\end{align}
and,  again  by the discrete Gronwall inequality \eq{9.2++} for a{}ny 
sufficiently small $\tau$, the estimate
\begin{align}
\|\overlinexitau\|_{L^\infty(I;L^s(\varOmega;\R^3))}^{}\le C\,.
\end{align}
Furthermore, the test of \eq{EUL-L-disc3}  
${\rm div}(|\nabla\bm\xi_{\etau}^k|^{ r_1-2}\nabla\bm\xi_{\etau}^k)$ with some
$ r_1>1$ needs to  handle  the convective term, which can be
performed as \eq{calulus-nonlin-rho-r} but for the vector-valued field, which
just replaces the tensorial product $\otimes$ by the tensorial product
$\boxtimes$ as the matrix $[(\nabla\bm\xi{\boxtimes}\nabla\bm\xi)]_{ij}
=\sum_{k=1}^3\Frac{\partial}{\partial x_i}\bm\xi_k\Frac{\partial}{\partial x_j}\bm\xi_k$.
Thus, like \eq{est-of-rho-disc} but here simpler, we obtain the inequality 
\begin{align}\nonumber
\Frac1{ r_1}\int_\varOmega\frac{|\nabla\bm\xi_{\etau}^k|^{ r_1}{-}|\nabla\bm\xi_{\etau}^{k-1}|^{ r_1}
\!\!}\tau\,\d\xx
\le\!\int_\varOmega\,&\Frac1{ r_1}({\rm div}\,\vv_{\etau}^k)|\nabla\bm\xi_{\etau}^k|^{ r_1}\!
-|\nabla\bm\xi_{\etau}^k|^{ r_1-2}(\nabla\bm\xi_{\etau}^k{\boxtimes}\nabla\bm\xi_{\etau}^k)\Colon\strain(\vv_{\etau}^k)\,\d\xx
\\[-.4em]&\qquad\quad
\le\Big(1{+}\Frac1{ r_1}\Big)\big\|\strain(\vv_{\etau}^k)\big\|_{L^\infty(\varOmega;\R^{3\times3})}^{}
\int_\varOmega|\nabla\bm\xi_{\etau}^k|^{ r_1}\,\d\xx\,\!\!
\end{align}
and by using again the discrete Gronwall inequality \eq{9.2++} relying on
\eq{Euler-small-est1}, for sufficiently small time steps $\tau$, we obtain
the estimate
\begin{align}\label{est-of-xi-disc}
\|\overlinexitau\|_{L^\infty(I;W^{1, r_1}(\varOmega;\R^3))}^{}\le C\,.
\end{align}

Eventually, we can perform similar tests also for \eq{EUL-L-disc2}.
First test is by $|\Fek|^{s-2}\Fek$. Using the calculus \eq{transport-rho}
for the convective term $(\vv_\etau^k\Cdot\nabla)\Fek$, we obtain
\begin{align}\nonumber
&\Frac1s\int_\varOmega\frac{|\Fek|^s{-}|\Fekk|^s\!\!}\tau\,\d\xx
\le\int_\varOmega\bigg(
\Big(\nabla\vv_\etau^k\!{+}\Frac1s{\rm div}\,\vv_\etau^k\Big)|\Fek|^s
-\Fek\mathscr{L}_\lambda^{}(\bm\xi_{\etau}^{k},\Fek)\Colon
(|\Fek|^{s-2}\Fek)\bigg)\,\d\xx
\\[-.3em]&\hspace{5em}
\le\Big(\big(1{+}\Frac1s\big)\|\nabla\vv_\etau^k\|_{L^\infty(\varOmega;\R^{3\times3})}^{}
+\|\mathscr{L}_\lambda^{}(\bm\xi_{\etau}^{k},\Fek)\|_{L^\infty(\varOmega;\R^{3\times3})}^{}\Big)
\int_\varOmega|\Fek|^s\,\d\xx\,.
\label{EUL-L-heter-Gronwall-E}\end{align}
Due to the truncation \eq{cutoff-stress} together with the
assumptions (\ref{EUL-L-ass}a,b), we can rely on that
$\sup_{X\in\varOmega,\ F\in\GL_3^+}|\mathscr{L}_\lambda^{}(X,F)|$ is finite;
note that the quadratic coercivity of $\zetap(X,\cdot)$
assumed in \eq{EUL-L-ass-zeta} ensures at most a linear
growth of $[\zetap(X,\cdot)^*]'$ which makes it bounded
on the truncated Mandel stress, so eventually it makes $|\mathscr{L}_\lambda^{}|$
bounded. Thus, from \eq{EUL-L-heter-Gronwall-E} by
the discrete Gronwall inequality for sufficiently small 
$\tau>0$ relying again on \eq{Euler-small-est1}, we obtain
\begin{align}
\|\overlineFetau\|_{L^\infty(I;L^s(\varOmega;\R^{3\times3}))}^{}\le C\,.
\end{align}
Furthermore, we perform the test of \eq{EUL-L-disc2} by
${\rm div}(|\nabla\Fek|^{r-2}\nabla\Fek)$. The convective term
$(\vv_{\eetau}^k\Cdot\nabla)\Fek$ can be treated as in
\eq{calulus-nonlin-rho-r}. Thus we get the inequality
\begin{align}\nonumber
\Frac1r&\!\int_\varOmega\frac{|\nabla\Fek|^r{-}|\nabla\Fekk|^r
\!\!}\tau\,\d\xx
\le\!\int_\varOmega\bigg(\Frac1r({\rm div}\,\vv_{\etau}^k)|\nabla\Fek|^r\!
-|\nabla\Fek|^{r-2}(\nabla\Fek{\boxtimes}\nabla\Fek)\Colon\strain(\vv_{\etau}^k)
\\[-.3em]\nonumber&\hspace{3.7em}+(\nabla\vv_{\etau}^k)|\nabla\Fek|^r
+(\nabla^2\vv_{\etau}^k)\Fek
-\nabla\big(\Fek\mathscr{L}_\lambda^{}(\bm\xi_{\etau}^{k},\Fek)\big)
\Vdots\big(|\nabla\Fek|^{r-2}\nabla\Fek\big)
\bigg)\,\d\xx
\\[-.3em]&\nonumber\le\Big(
\frac{2r{+}1}r\big\|\nabla\vv_{\etau}^k\big\|_{L^\infty(\varOmega;\R^{3\times3})}^{}\!
+\!\!\sup_{X\in\varOmega,\,F\in\R^{3\times3}}\big|\mathscr{L}_\lambda^{}(X,F)\big|
+\big|F[\mathscr{L}_\lambda^{}]_F'(X,F)\big|\Big)
\!\int_\varOmega\!\big|\nabla\Fek|^r\d\xx
\\[-.3em]&\nonumber\hspace{9.7em}
+\|\Fek\|_{L^{pr/(p-r)}}^{}\|\nabla^2\vv_{\etau}^k\|_{L^p(\varOmega;\R^{3\times3\times3})}^{}
\|\nabla\Fek\|_{L^r(\varOmega;\R^{3\times3\times3})}^{r-1}
\\[-.3em]&\hspace{9.7em}
-\!\int_\varOmega\Big(\Fek\big[\mathscr{L}_\lambda^{}\big]_X'(\bm\xi_{\etau}^{k},\Fek)\big)
\nabla\bm\xi_{\etau}^{k}\Big)\Vdots\big(|\nabla\Fek|^{r-2}\nabla\Fek\big)\,\d\xx
\,;\!\!
\label{EUL-L-heter-Gronwall+}\end{align}
actually, here the tensorial product $\boxtimes$ means
the matrix $[\nabla\FF{\boxtimes}\nabla\FF]_{ij}=\sum_{k,l=1}^3
\Frac{\partial}{\partial x_i}\FF_{kl}\Frac{\partial}{\partial x_j}\FF_{kl}$.
Estimating of the last term relies on the boundedness of
$\sup_{X\in\varOmega,\,F\in\R^{3\times3}}|F\,[\mathscr{L}_\lambda^{}]_X'(X,F)|$  for
a fixed $\lambda$ and on \eq{est-of-xi-disc} with $r<p$ 
since $\Fezero\in W^{1, r}(\varOmega;\R^{3\times3})$. Thus this
 last  term in
\eq{EUL-L-heter-Gronwall+} can be estimated up to a multiplicative
constant by $1+\|\nabla\Fek\|_{L^{r}(\varOmega;\R^{3\times3\times3})}^{r}$.
Then, using the already obtained estimates \eq{Euler-small-est1} and
\eq{est-of-xi-disc},
we can use
again the discrete Gronwall inequality \eq{9.2++} so that, for sufficiently
small time steps $\tau$, we obtain the estimate
\begin{align}\label{est-of-E-disc}
\|\overlineFetau\|_{L^\infty(I;W^{1,r}(\varOmega;\R^{3\times3}))}^{}\le C\,.
\end{align}

Furthermore, by comparison from \eq{EUL-L-visco-ED+discr} and
the already obtained estimates, we still obtain 
\begin{subequations}\label{EUL-L-est-dt}
\begin{align}\label{EUL-L-est-dt-rho}
&\Big\|\pdt{\varrho_{\etau}}\Big\|_{L^p(I;L^r(\varOmega))}^{}\le C\,,
\\&\Big\|\pdt{\pp_{\etau}}\Big\|_{L^{p'}(I;W^{2,p}(\varOmega;\R^3)^*)}^{}\le C\,,
\\&\Big\|\pdt{\Fe_{,\tau}}\Big\|_{L^p(I;L^r(\varOmega;\R^{3\times3}))}^{}\le C\,,\
\text{ and}
\\&\Big\|\pdt{\bm\xi_\tau}\Big\|_{L^p(I;L^{ r_1}(\varOmega;\R^3))}^{}\le C\ \ \text{ for any }\ 1\le r_1<\infty\,.
\end{align}\end{subequations}
These estimates hold also for the piecewise-constant interpolants
if $L^p(I;\cdot)$ or $L^{p'}(I;\cdot)$ are replaced by the space of
vector-valued measures.

\medskip{\it Step 3: Existence of a solution to \eq{EUL-L-disc}}.
The rigorous proof of existence of a solution $(\varrho_\tau^k,\vv_\tau^k)$ of
the subsystem (\ref{EUL-L-disc}a,b) is slightly delicate.
Convincing arguments  can  rely on a suitable regularization
towards  a  quasilinear elliptic system  as specified in
\eq{EUL-L-disc-heter} below, 
for which standard theory can be used, and
the{n}  a subsequent limit passage.

Due to \eq{Euler-small-est1} with $p>3$, we know, in particular, that
\begin{align}\label{EUL-L-est-div-v}
\max_{1\le k\le T/\tau}\Big\|\frac{\pp_\tau^k}{\sqrt{\varrho_\tau^k}}\Big\|_{L^2(\varOmega;\R^3)}\le C\ \ \text{ and }\ \ 
\tau\sum_{k=1}^{T/\tau}\|{\rm div}\,\vv_\tau^k\|_{L^\infty(\varOmega)}^2\le C\,.
\end{align}

For a given $(\varrho_{\tau}^{k-1},\pp_{\tau}^{k-1},\Fekk,\bm\xi_{\tau}^{k-1})\in
 W^{1,r}(\varOmega;\R\times\R^3\times\R^{3\times3}\times\R^3)$
with $\varrho_{\tau}^{k-1}>0$ a.e.,
we seek a solution $(\varrho_{\EPS\DELTA\tau}^k,\vv_{\EPS\DELTA\tau}^k)$ and thus
also $\pp_{\EPS\DELTA\tau}^k$ of the regularized quasilinear elliptic system 
\begin{subequations}\label{EUL-L-disc-heter}
\begin{align}\label{EUL-L-0disc-heter}
&\!\!\frac{\varrho_{\EPS\DELTA\tau}^k{-}\varrho_\tau^{k-1}\!\!}\tau\,=
{\rm div}\big(\DELTA|\nabla\varrho_{\EPS\DELTA\tau}^k|^{r-2}
\nabla\varrho_{\EPS\DELTA\tau}^k\!-\pp_{\EPS\DELTA\tau}^k\big)\ \ 
\text{ with }\,\pp_{\EPS\DELTA\tau}^k=\varrho_{\EPS\DELTA\tau}^k\vv_{\EPS\DELTA\tau}^k
\ \ \text{ and}
\\[-.2em]&\nonumber
\!\!\frac{\pp_{\EPS\DELTA\tau}^k{-}\pp_\tau^{k-1}\!\!}\tau\,=
{\rm div}\Big(\mathscr{T}_\lambda(\bm\xi_{\tau}^{k-1},\Fekk)
{+}\DD_{\EPS\DELTA\tau}^k{-}\pp_{\EPS\DELTA\tau}^k{\otimes}\vv_{\EPS\DELTA\tau}^k\Big)
+\varrho_{\EPS\DELTA\tau}^k\GRAVITY_{\tau}^k
\\[-.3em]&\nonumber\hspace*{13em}
-\EPS|\vv_{\EPS\DELTA\tau}^k|^{p-2}\vv_{\EPS\DELTA\tau}^k-\DELTA
|\nabla\varrho_{\EPS\DELTA\tau}^k|^{r-2}(\nabla\vv_{\EPS\DELTA\tau}^k)
\nabla\varrho_{\EPS\DELTA\tau}^k
\\[-.2em]&
\hspace*{1em}\text{ with }\ \ \DD_{\EPS\DELTA\tau}^k=\bbD\strain(\vv_{\EPS\DELTA\tau}^k)
-{\rm div}\,\mathfrak{H}_{\EPS\DELTA\tau}^k\ \ \text{ with }\ \
\mathfrak{H}_{\EPS\DELTA\tau}^k=\HYPER\big|\nabla^2\vv_{\EPS\DELTA\tau}^k\big|^{p-2}\nabla^2\vv_{\EPS\DELTA\tau}^k\!
\label{EUL-L-1disc-heter}
\end{align}\end{subequations}
considered with $r>3$  and with the additional boundary condition
$\nn\Cdot\nabla\varrho_{\EPS\DELTA\tau}^k=0$; for the $\DELTA$-terms in the
case $r=2$ (not used here) we refer to \cite{Zato12ASCN}.
Like \eq{EUL-L-basic-engr-balance-disc}, we obtain the inequality
\begin{align}\nonumber
\!\!\!\!\!&\int_\varOmega\!\frac{|\pp_{\EPS\DELTA\tau}^k|^2\!\!\!}{2\varrho_{\EPS\DELTA\tau}^k\!\!\!}
\,\d\xx
+\tau\!\int_\varOmega\!
\bbD\strain(\vv_{\EPS\DELTA\tau}^k)\Colon\strain(\vv_{\EPS\DELTA\tau}^k)
+\HYPER|\nabla^2\vv_{\EPS\DELTA\tau}^k|^p+\EPS|\vv_{\EPS\DELTA\tau}^k|^p
\,\d\xx
\\[-.5em]&\hspace{7em}
\le\int_\varOmega\!\frac{|\pp_{\tau}^{k-1}|^2\!\!}{2\varrho_{\tau}^{k-1}\!\!}\,\d\xx
+\tau\!
\int_\varOmega\!\varrho_{\EPS\DELTA\tau}^k\GRAVITY_{\tau}^k\Cdot\vv_{\EPS\DELTA\tau}^k
-\mathscr{T}_\lambda^{}(\bm\xi_{\tau}^{k-1},\Fekk)\Colon\strain(\vv_{\EPS\DELTA\tau}^k)
\,\d\xx\,;\!\!
\label{EUL-L-basic-engr-balance-disc+}\end{align}
here we used cancellation of the two $\DELTA$-regularizing terms as in
\cite{Roub25TDVE,Zato12ASCN}. Moreover, testing \eq{EUL-L-0disc-heter}
by $\varrho_{\EPS\DELTA\tau}^k$, we obtain
\begin{align}\nonumber
\int_\varOmega\!\DELTA|\nabla\varrho_{\EPS\DELTA\tau}^k|^r+\Frac1\tau(\varrho_{\EPS\DELTA\tau}^k)^2\,\d\xx
&\le\int_\varOmega\!\Frac1\tau\varrho_{\EPS\DELTA\tau}^k\varrho_\tau^{k-1}
+\varrho_{\EPS\DELTA\tau}^k\vv_{\EPS\DELTA\tau}^k\Cdot\nabla\varrho_{\EPS\DELTA\tau}^k\,\d\xx
\\&=\int_\varOmega\!\Frac1\tau\varrho_{\EPS\DELTA\tau}^k\varrho_\tau^{k-1}
-\Frac12({\rm div}\,\vv_{\EPS\DELTA\tau}^k)\,(\varrho_{\EPS\DELTA\tau}^k)^2\,\d\xx\,,
\label{EUL-L-test-contiunity-reg}\end{align}
where we used the calculus \eq{transport-rho} for $s=2$. 
Applying the Young inequality to \eq{EUL-L-basic-engr-balance-disc+}
and \eq{EUL-L-test-contiunity-reg} shows a-priori estimates for
$(\varrho_{\EPS\DELTA\tau}^k,\vv_{\EPS\DELTA\tau}^k)$ in $W^{1,r}(\varOmega)\times W^{2,p}(\varOmega;\R^3)$,
which yields existence of this solution by standard methods for quasilinear
elliptic systems. In particular, we have also the estimate
\begin{align}\label{EUL-L-est-v-regularized}
\|\vv_{\EPS\DELTA\tau}^k\|_{L^p(\varOmega;\R^3)}^{}=C/{\sqrt[p]\EPS}\,.
\end{align}

Moreover, testing \eq{EUL-L-0disc-heter} by 
${\rm div}(|\nabla\varrho_{\EPS\DELTA\tau}^k|^{r-2}\nabla\varrho_{\EPS\DELTA\tau}^k)$
 which eventually can be shown to belong to $L^2(\varOmega)$ (after
limiting a Galerkin approximation to be truly rigorous at this point) 
 gives the estimate
$\|\nabla\varrho_{\EPS\DELTA\tau}^k\|_{L^r(\varOmega;\R^3)}^{}\le C$ uniformly
in $\DELTA$. This allows for the limit passage with $\DELTA\to0$, which gives
a solution
$(\varrho_{\EPS\tau}^k,\vv_{\EPS\tau}^k)\in W^{1,r}(\varOmega)\times W^{2,p}(\varOmega;\R^3)$
of the quasilinear elliptic system
\begin{subequations}\label{EUL-L-disc-heter+}
\begin{align}\label{EUL-L-0disc-heter+}
&\!\!\frac{\varrho_{\EPS\tau}^k{-}\varrho_\tau^{k-1}\!\!}\tau\,=
-{\rm div}\,\pp_{\EPS\tau}^k\ \ 
\text{ with }\,\pp_{\EPS\tau}^k=\varrho_{\EPS\tau}^k\vv_{\EPS\tau}^k\,,\!
\\[-.2em]&\nonumber
\!\!\frac{\pp_{\EPS\tau}^k{-}\pp_\tau^{k-1}\!\!}\tau\,=
{\rm div}\Big(\mathscr{T}_\lambda(\bm\xi_{\tau}^{k-1},\Fekk)
{+}\DD_{\EPS\tau}^k{-}\pp_{\EPS\tau}^k{\otimes}\vv_{\EPS\tau}^k\Big)
+\varrho_{\EPS\tau}^k\GRAVITY_{\tau}^k
-\EPS|\vv_{\EPS\tau}^k|^{p-2}\vv_{\EPS\tau}^k\,,
\\[-.3em]&
\hspace*{4em}\text{ where }\ \DD_{\EPS\tau}^k=\bbD\strain(\vv_{\EPS\tau}^k)
-{\rm div}\,\mathfrak{H}_{\EPS\tau}^k\ \ \text{ with }\ \
\mathfrak{H}_{\EPS\tau}^k=\HYPER\big|\nabla^2\vv_{\EPS\tau}^k\big|^{p-2}\nabla^2\vv_{\EPS\tau}^k\,.\!
\label{EUL-L-1disc-heter+}
\end{align}\end{subequations}
For the obtained $\vv_{\EPS\tau}^k\in W^{2,p}(\varOmega;\R^3)$, the mass continuity
equation \eq{EUL-L-0disc-heter+} has a unique solution
$\varrho_{\EPS\tau}^k\in W^{1,r}(\varOmega)$. By testing \eq{EUL-L-0disc-heter+}
by the negative part of $\varrho_{\EPS\tau}^k$, we can see that 
$\varrho_{\EPS\tau}^k\ge0$. Moreover, even $\varrho_{\EPS\tau}^k>0$ a.e.,  
which can be seen by a contradiction argument:
assuming $\varrho_\tau^{k-1}>0$ on $\barOmega$ and the minimum of a
(momentarily smooth) solution $\varrho_{\EPS\tau}^k$ 
is attained at some $\xx\in\varOmega$ and $\varrho_{\EPS\tau}^k(\xx)=0$,
so that $\varrho_{\EPS\tau}^k(\xx)-\varrho_\tau^{k-1}(\xx)<0$ and
$\nabla\varrho_{\EPS\tau}^k(\xx)=\bm0$. Then
$\varrho_{\EPS\tau}^k(\xx)-\varrho_\tau^{k-1}(\xx)=
\tau\vv_{\EPS\tau}^k(\xx)\Cdot\nabla\varrho_{\EPS\tau}^k(\xx)
-\tau({\rm div}\,\vv_{\EPS\tau}^k(\xx))\varrho_{\EPS\tau}^k(\xx)\ge0$, so that we
obtain the desired contradiction showing that $\varrho_{\EPS\tau}^k(\xx)>0$.

This a.e.\ positivity allows for a test by
$\sigma_{\EPS\tau}^k:=1/\varrho_{\EPS\tau}^k$. Using the convexity of
$\varrho\mapsto1/\varrho$ on $(0,+\infty)$ for the convective time differences
and \eq{EUL-L-0disc-heter+}, we obtain a.e.\ on $\varOmega$ the inequality
\begin{align}\nonumber
&\frac{\sigma_{\EPS\tau}^k-\sigma_{\tau}^{k-1}\!\!\!}\tau
+\vv_{\EPS\tau}^k\Cdot\nabla\sigma_{\EPS\tau}^k
\ =\frac1\tau\Big(\frac1{\varrho_{\EPS\tau}^k}-\frac1{\varrho_{\tau}^{k-1}}\Big)
+\vv_{\EPS\tau}^k\Cdot\nabla\frac1{\varrho_{\EPS\tau}^k}
\\&\qquad
\le-\frac1{(\varrho_{\EPS\tau}^k)^2\!}\Big(\frac{\varrho_{\EPS\tau}^k-\varrho_{\tau}^{k-1}\!\!\!}\tau
+\vv_{\EPS\tau}^k\Cdot\nabla\varrho_{\EPS\tau}^k\Big)
=\frac1{\varrho_{\EPS\tau}^k\!}{\rm div}\,\vv_{\EPS\tau}^k=({\rm div}\,\vv_{\EPS\tau}^k)\,\sigma_{\EPS\tau}^k\,,
\label{ineq-for-sparsity}\end{align}
where $\sigma_{\tau}^{k-1}:=1/\varrho_\tau^{k-1}$. Testing  \eq{ineq-for-sparsity}
by $|\sigma_{\EPS\tau}^k|^{s-2}\sigma_{\EPS\tau}^k$ for some $s>1$ 
and using the calculus \eq{transport-rho} now modified, this test gives the inequality
\begin{align}\nonumber
\frac1s\int_\varOmega\frac{(\sigma_{\EPS\tau}^k)^s
-(\sigma_{\tau}^{k-1})^s}\tau\,\d\xx
&\le\int_\varOmega\Big(1{+}\frac1s\Big)({\rm div}\,\vv_{\EPS\tau}^k)
(\sigma_{\EPS\tau}^k)^s\,\d\xx
\\&
\le\Big(1{+}\frac1s\Big)\|{\rm div}\,\vv_{\EPS\tau}^k\|_{L^\infty(\varOmega)}^{}
\int_\varOmega(\sigma_{\EPS\tau}^k)^{s}\,\d\xx\,.
\label{est-of-sparsity}\end{align}
For a{}ny  sufficiently small $\tau>0$, say for
$\tau\le1/(\sqrt2\max_{1\le k\le T/\tau}
\|{\rm div}\,\vv_{\EPS\tau}^k\|_{L^\infty(\varOmega)}^{})$, we obtain a uniform
bound for $\|\sigma_{\EPS\tau}^k\|_{L^s(\varOmega)}^{}$, cf.\ the discrete
Gronwall inequality \eq{9.2++}. Here note that $\max_{1\le k\le T/\tau}
\|{\rm div}\,\vv_{\EPS\tau}^k\|_{L^\infty(\varOmega)}^{}\le\sqrt{C/\tau}$ with
$C$ from the second estimate in \eq{EUL-L-est-div-v} so that 
choosing $\tau\le 1/(2C)$ is sufficient. Having
$\|\sigma_{\EPS\tau}^k\|_{L^s(\varOmega)}^{}$ bounded and using the first estimate
in \eq{EUL-L-est-div-v} while realizing that $\vv_{\EPS\tau}^k=
\sqrt{\sigma_{\EPS\tau}^k}\,\pp_{\EPS\tau}^k/\sqrt{\varrho_{\EPS\tau}^k}$,
we obtain a bound 
\begin{align}
\|\vv_{\EPS\tau}^k\|_{L^{2s/(s+1)}(\varOmega;\R^3)}^{}
\le C_s\big\|\sqrt{\sigma_{\EPS\tau}^k}\big\|_{L^{2s}(\varOmega)}^{}
\Big\|\frac{\pp_{\EPS\tau}^k}{\sqrt{\varrho_{\EPS\tau}^k}}\Big\|_{L^2(\varOmega;\R^3)}\,.
\label{EUL-L-est-v-2-eps}\end{align}
Having $\{\vv_{\EPS\tau}^k\}_{\EPS>0}$ bounded in $L^{2s/(s+1)}(\varOmega;\R^3)$,
we have also the bound in $W^{2,p}(\varOmega;\R^3)$.
We can pass to the limit with $\EPS\to0$. The regularizing term
$\EPS|\vv_{\EPS\tau}^k|^{p-2}\vv_{\EPS\tau}^k$ in \eq{EUL-L-1disc-heter+},
whose norm in $L^{p'}(\varOmega;\R^3)$ is $\mathscr{O}(\EPS^{1/(p-1)})$, 
obviously vanishes for $\EPS\to0$. Thus the existence of a solution
$(\varrho_\tau^k,\vv_\tau^k)\in W^{1,r}(\varOmega)\times W^{2,p}(\varOmega;\R^3)$ of
(\ref{EUL-L-disc}a,b) is proved. In fact, choosing $s=2$ leading
to the boundedness of $\{\vv_{\EPS\tau}^k\}_{\EPS>0}$ in $L^{4/3}(\varOmega;\R^3)$
would be sufficient for it. Moreover, from \eq{EUL-L-est-v-2-eps} with
an arbitrarily large $s>1$, we obtained an estimate
\begin{align}
\|\overlinevvtau\|_{L^\infty(I;L^a(\varOmega;\R^3))}^{}\le C_a\ \ \text{ with
an arbitrary }\  1\le a<2\,.
\label{EUL-L-est-v-2-eps+}\end{align}
Together with \eq{Euler-small-est1}, we thus obtain the bound for
$\overlinevvtau$ in $L^a(I;W^{2,p}(\varOmega;\R^3))$ with $a<2$,
in particular 
\begin{align}
\|\overlinevvtau\|_{L^{p}(I;W^{1,\infty}(\varOmega;\R^3))}^{}\le C\,.
\label{EUL-L-est-v-2-eps++}\end{align}
Moreover, \eq{EUL-L-est-v-2-eps+} together with
\eq{Euler-small-est1} and \eq{est-of-rho-disc} allow us also to augment 
\eq{Euler-small-est3+++} by an estimate for $\nabla\overlineppetau=
\nabla(\overlinerhoetau\overlinevvtau)=\overlinerhoetau\nabla\overlinevvtau
+\nabla\overlinerhoetau{\otimes}\overlinevvtau$, namely
\begin{align}
\|\overlineppetau\|_{L^{p}(I;W^{1,r}(\varOmega;\R^3))}^{}\le C\,.
\label{EUL-L-est-p+}\end{align}

Noteworthy, \eq{ineq-for-sparsity} is at disposal only as an inequality,
therefore, unlike \eq{EUL-L-disc0}, we cannot perform the test by
${\rm div}(|\nabla\sigma_{\EPS\tau}^k|^{r-2}\nabla\sigma_{\EPS\tau}^k)$ to
obtain an estimate of $\nabla\sigma_{\EPS\tau}^k$.  On the other hand,
we can read such an estimate from 
$\nabla\sigma_{\EPS\tau}^k=\nabla(1/\varrho_{\EPS\tau}^k)=
-(\sigma_{\EPS\tau}^k)^2\nabla\varrho_{\EPS\tau}^k\in W^{1,rs/(2r+s)}(\varOmega;\R^3)$.
In the limit, realizing also that $r>3$ and that $s\ge1$ can be considered
arbitrarily large, we obtain the bound of $\osetau$ in
$L^\infty(I{\times}\varOmega)$. From \eq{Euler-small-est1}, we can then read
the bound for $\overlineppetau$ in $L^\infty(I;L^2(\varOmega;\R^3))$ and thus
also for $\overlinevvtau=\osetau\overlineppetau$ in
$L^\infty(I;L^2(\varOmega;\R^3))$.

\medskip{\it Step 4: Limit passage for $\tau\to0$}.
By the Banach selection principle, we obtain a subsequence converging
weakly* with respect to the topologies indicated in \eq{Euler-small-est}, 
\eq{est-of-xi-disc}, \eq{est-of-E-disc}, \eq{EUL-L-est-dt}, and
\eq{EUL-L-est-p+}
to some limit
$(\varrho_\EEps,\pp_\EEps,\vv_\EEps,\Fe_\EEps,\bm\xi_\EEps)$. Specifically,
\begin{subequations}\label{Euler-small-converge}
\begin{align}
&&&\!\!\overlinerhoetau\to\varrho_\EEps&&\text{weakly* in $\ L^\infty(I;W^{1,r}(\varOmega))$\,,}
\\&&&\!\!\varrho_\etau\to\varrho_\EEps&&\text{weakly* in $\ L^\infty(I;W^{1,r}(\varOmega))\,\cap\, 
W^{1,p}(I;L^r(\varOmega))$\,,}
\\&&&\label{Euler-small-converge-bar-p}
\!\!\overlineppetau\to\pp_\EEps&&\text{weakly\ \;in $\
L^p(I;W^{1,r}(\varOmega;\R^3))$\,,}
\\&&&\label{Euler-small-converge-p}
\!\!\pp_\etau\to\pp_\EEps&&\text{weakly\ \;in $\
L^p(I;W^{1,r}(\varOmega;\R^3))\,\cap\,W^{1,p'}(I;W^{2,p}(\varOmega;\R^3)^*)$\,,}
\\&&&\label{Euler-small-converge-bar-v}
\!\!\overlinevvtau\to\vv_\EEps&&\text{weakly\ \;in $\ L^{2}(I;W^{2,p}(\varOmega;\R^3))$\,,}
\\[-.1em]
&&&\!\!\overlineFetau\to\Fe_\EEps\ \text{ and }\
\underlineFetau\to\Fe_\EEps\!\!\!\hspace{-7em}&&\hspace{7em}\text{weakly* in $\ 
L^\infty(I;W^{1,r}(\varOmega;\R^{3\times3}))$\,,}\!\!
\\
&&&\!\!\FF_{{\rm e},\tau}\to\Fe_\EEps\!\!\!&&\text{weakly*\ \;in $\ 
L^\infty(I;W^{1,r}(\varOmega;\R^{3\times3}))\,\cap\,W^{1,p}(I;L^r(\varOmega;\R^{3\times3}))$\,,}\!\!\!&&
\\[-.1em]
&&&\!\!\overlinexitau\to\bm\xi_\EEps\ \text{ and }\
\underlinexitau\to\bm\xi_\EEps\!\!\!\hspace{-6em}&&\hspace{5em}\text{weakly* in $\ 
L^\infty(I;W^{1, r_1}(\varOmega;\R^{3\times3}))$
\,,}\!\!
\\&&&\!\!
\bm\xi_\etau\to\bm\xi_\EEps\!\!\!&&\text{weakly\ \;in $\ 
L^\infty(I;W^{1, r_1}(\varOmega;\R^3))\,\cap\,W^{1,p}(I;L^{ r_1}(\varOmega;\R^3))$\,;}
\end{align}\end{subequations}
here the exponent $r$ in (\ref{Euler-small-converge}a-d,f,g) refers to 
$3< r<p$ from \eq{EUL-L-ass-IC} while in
(\ref{Euler-small-converge}h,i)
the exponent  $r_1<\infty$ referring to \eq{est-of-xi-disc} 
can be arbitrarily large.
Notably, the limit of $\overlinerhoetau$ and $\varrho_\etau$ is indeed the same
due to the control of $\pdt{}\varrho_\etau$ in \eq{EUL-L-est-dt-rho}; cf.\
\cite[Sect.8.2]{Roub13NPDE}. The same is true also for $\overlineppetau$ and
$\pp_\etau$, and for $\overlineEetau$  and $\underlineFetau$  
and $\FF_{{\rm e},\tau}$,  and similarly for $\bm\xi$'s, too.

By the compact embedding $W^{1,r}(\varOmega)\subset C(\barOmega)$ for $r>3$
and the (generalized) Aubin-Lions theorem, cf.\ \cite[Corollary~7.9]{Roub13NPDE},
we have also, for arbitrary $1\le a<\infty$, 
\begin{subequations}\label{Euler-small-converge-strong}
\begin{align}\label{Euler-small-converge-strong-rho}
&&&\!\!\overlinerhoetau\to\varrho_\EEps\!\!\!&&\text{strongly in }\
L^a(I;C(\barOmega))\,,&&
\\\label{Euler-small-converge-strong-p-}
&&&\!\!\overlineppetau\to\pp_\EEps\!\!\!\!\!\!&&\text{strongly in }\ L^p(I;C(\barOmega;\R^3))\,,
\\\label{Euler-small-converge-strong-E}
&&&\!\!\overlineEetau\!\to\Fe_\EEps\ \text{ and }\
\underlineFetau\!\to\Fe_\EEps\!\!\!\hspace{-8em}&&\hspace{6em}
\text{strongly in }\
L^a(I;C(\barOmega;\R^{3\times3}))\,, \text{ and }
\\\label{Euler-small-converge-strong-xi}&&&
\!\!\overlinexitau\!\to\bm\xi_\EEps\ \ \ \ \text{ and }\ \ 
\underlinexitau\to\bm\xi_\EEps\!\!\!\hspace{-8em}&&\hspace{6em}
\text{strongly in }\ L^a(I;C(\barOmega;\R^{3}))\,.
\intertext{Moreover, using the Arzel\`a-Ascoli-type theorem, cf.\
\cite[Lemma~7.10]{Roub13NPDE}, we have also}
\label{Euler-small-converge-strong-rho+}
&&&\!\!\varrho_\etau\to\varrho_\EEps\!\!\!&&\text{strongly in }\ C(I{\times}\barOmega)\,.
\end{align}\end{subequations}

By \eq{Euler-small-converge-strong-E} and \eq{Euler-small-converge-strong-xi}
and by the continuity and boundedness of
$\mathscr{T}_\lambda,\mathscr{L}_\lambda:\varOmega\times\R^{3\times3}\to\R^{3\times3}$,
we can pass to the limit in the nonlinear terms:
\begin{subequations}\label{Euler-small-converge-strong+}
\begin{align}\label{Euler-small-converge-strong-T}
&\mathscr{T}_\lambda(\underlinexitau,\underlineFetau)
\to\mathscr{T}_\lambda(\bm\xi_\EEps,\Fe_\EEps)
&&\hspace{-0em}\text{strongly\ in }\, 
L^a(I{\times}\varOmega;\R^{3\times3})\,,
\\&\mathscr{L}_\lambda(\overlinexitau,\overlineEetau)
\to\mathscr{L}_\lambda(\bm\xi_\EEps,\Fe_\EEps)
&&\hspace{-0em}\text{strongly\ in }\,
L^a(I{\times}\varOmega;\R^{3\times3})\ \text{ for any }1\le a<\infty\,.
\label{Euler-small-converge-strong+L}
\end{align}\end{subequations}

This already allows for passing to the limit in
\eq{EUL-L-visco-ED+discr}. While the limit passage
in \eq{EUL-L-visco-ED+0discr} is simple
due to (\ref{Euler-small-converge}b,c), the limit passage
in the quasilinear momentum equation \eq{EUL-L-visco-ED+1discr}
is a bit more technical. Using the monotonicity of the operator
$\vv\mapsto{\rm div}({\rm div}(\HYPER|\nabla^2\vv|^{p-2}\nabla^2\vv)
-\bbD\strain(\vv))$  understood, of course, in the
weak formulation together  with the boundary conditions
\eq{Euler-small-BC-hyper} and using \eq{EUL-L-visco-ED+1discr} tested
by $\ovetau{-}\wt\vv$, we obtain
\begin{align}\nonumber
\nonumber\hspace*{0em}
0&\le\!\int_0^T\!\!\!\int_\varOmega\!\bigg(\bbD\strain(\ovetau{-}\wt\vv)
\Colon\strain(\ovetau{-}\wt\vv)
 +\HYPER\Big(|\nabla^2\ovetau|^{p-2}\nabla^2\ovetau
-|\nabla^2\wt\vv|^{p-2}\nabla^2\wt\vv\Big)\Vdots\nabla^2(\ovetau{-}\wt\vv)
\bigg)\,\d\xx\d t
 \\[-.0em]&=\nonumber
 \int_0^T\!\!\!\int_\varOmega\bigg(
 \Big(\oretau\overline\GRAVITY_{\tau}-\pdt{\pp_\etau}\Big)
 \Cdot(\ovetau{-}\wt\vv)
-\Big(\mathscr{T}_\lambda(\underlinexitau,\underlineFetau)
{-}\opetau{\otimes}\ovetau\Big)\Colon\strain(\ovetau{-}\wt\vv)
\\[-.7em]&\nonumber\hspace{12.8em}
-\bbD\strain(\wt\vv)\Colon\strain(\ovetau{-}\wt\vv)
-\HYPER|\nabla^2\wt\vv|^{p-2}\nabla^2\wt\vv\Vdots\nabla^2(\ovetau{-}\wt\vv)
\bigg)\,\d\xx\d t
\\[-.5em]&\le\nonumber
 \int_0^T\!\!\!\int_\varOmega\bigg(
\oretau\overline\GRAVITY_{\tau}\Cdot(\ovetau{-}\wt\vv)
 +\pdt{\pp_\etau}\Cdot\wt\vv
-(\opetau{\otimes}\ovetau)\Colon\strain(\wt\vv)
-\big(\bbD\strain(\wt\vv){+}
\mathscr{T}_\lambda(\underlinexitau,\underlineFetau)\big)
\Colon\strain(\ovetau{-}\wt\vv)
\\[-.4em]&\hspace{8em}
-\HYPER|\nabla^2\wt\vv|^{p-2}\nabla^2\wt\vv\Vdots\nabla^2(\ovetau{-}\wt\vv)
\bigg)\,\d\xx\d t+
\int_\varOmega\frac{|\pp_0|^2}{2\varrho_0}-\frac{|\pp_\etau(T)|^2}{2\varrho_\etau(T)}\,\d\xx
\label{EUL-L-strong+}\end{align}
for any $\wt\vv\in L^p(I;W^{2,p}(\varOmega;\R^3))$. 
The last inequality in \eq{EUL-L-strong+} has again exploited the
convexity of the kinetic energy $(\pp,\varrho)\mapsto\frac12|\pp|^2/\varrho$
in the calculus:
\begin{align}\nonumber
\int_\varOmega&\frac{|\pp_\etau(T)|^2\!}{2\varrho_\etau(T)}
-\frac{|\pp_0|^2}{2\varrho_0}\,\d\xx
\le\int_0^T\!\!\!\int_\varOmega\!\pdt{\pp_\etau\!}\Cdot\ovetau
-\frac{|\ovetau|^2}2\pdt{\varrho_{\etau}}\,\d\xx\d t
\\&
{\buildrel{\scriptsize\eq{EUL-L-visco-ED+0discr}}\over{=}}
\int_0^T\!\!\!\int_\varOmega\!\pdt{\pp_\etau\!}\Cdot\ovetau
+\frac{|\ovetau|^2\!\!}2\:{\rm div}\,\opetau\,\d\xx\d t
=
\int_0^T\!\!\!\int_\varOmega\!\pdt{\pp_\etau\!}\Cdot\ovetau
+\ovetau\Cdot{\rm div}\big(\opetau{\otimes}\ovetau\big)
\,\d\xx\d t\,,
\label{EUL-L-strong-calc}
\end{align}
where, for the last equality, we used the Green formula for the calculus 
\begin{align}\nonumber
\!\!\!\!\!\int_\varOmega\varrho(\vv\Cdot\nabla)\vv\Cdot\vv\,\d\xx
&=\int_\varGamma\varrho|\vv|^2\!\!\!\!\lineunder{\!\!\!\!\vv\Cdot\nn\!\!\!\!}{$=0$}\!\!\!\!\d S
-\!\!\int_\varOmega\!\vv\Cdot{\rm div}(\varrho\vv{\otimes}\vv)\,\d\xx
\\[-.5em]&
=-\!\!\int_\varOmega\varrho(\vv\Cdot\nabla)\vv\Cdot\vv+{\rm div}(\varrho\vv)|\vv|^2\,\d\xx
=-\int_\varOmega\!{\rm div}(\varrho\vv)\frac{|\vv|^2}2\,\d\xx\,.
\label{calculus-convective}\end{align}
Now we want to pass to the limit in \eq{EUL-L-strong+} or, more
precisely, to estimate the limit superior from above. For this, we again
use the kinetic-energy convexity, which causes the weak lower semicontinuity
of $(\varrho,\pp)\mapsto\int_\varOmega|\pp|^2/\varrho\,\d\xx$ as a convex
functional
$\{\rho\,{\in}\,L^1(\varOmega);\,\rho>0\}\times L^1(\varOmega;\R^3)
\to[0,+\infty]$.
Here we rely also on that $|\pp_\etau(T)|^2/\varrho_\etau(T)$ is bounded in
$L^1(\varOmega)$ due to the former estimate in \eq{Euler-small-est1} and on
that $\varrho_\etau(T)\to\varrho_\EEps(T)$ even strongly in $C(\barOmega)$ due
to \eq{Euler-small-converge-strong-rho+}, and on that $\pp_\etau(T)$
converges weakly* in $C(\barOmega)^*$, i.e.\ as measures on $\barOmega$ due
to \eq{Euler-small-est3+++} to its limit which is $\pp_\EEps(T)$ because
simultaneously $\pp_\etau(T)\to\pp_\EEps(T)$ weakly in $W^{2,p}(\varOmega;\R^3)^*$
due to \eq{Euler-small-converge-p}.
For the term $(\opetau{\otimes}\ovetau)\Colon\strain(\wt\vv)$,
we use \eq{Euler-small-converge-bar-v} with 
\eq{Euler-small-converge-strong-p-}, which implies in particular that
$\opetau{\otimes}\ovetau\to\pp{\otimes}\vv$ weakly in
$L^a(I{\times}\varOmega;\R^{3\times3})$ for any $1<a<2p/(p{+}2)$.
All this allows us to estimate the limit superior of
\eq{EUL-L-strong+} from above:
\begin{align}
\nonumber
0&\le\int_0^T\!\!\bigg(\bigg\langle\pdt{\pp_\EEps},\wt\vv\bigg\rangle
+\!\int_\varOmega\Big(\varrho_\EEps\GRAVITY_\EEps\Cdot(\vv_\EEps{-}\wt\vv)
-(\pp_\EEps{\otimes}\vv_\EEps)\Colon\strain(\wt\vv)
  -\big(\mathscr{T}_\lambda(\bm\xi,\Fe){+}\bbD\strain(\wt\vv)\big)
  \Colon\strain(\vv_\EEps{-}\wt\vv)
\\[-.7em]&\nonumber\hspace{9em}
-\HYPER|\nabla^2\wt\vv|^{p-2}\nabla^2\wt\vv\Vdots\nabla^2(\vv_\EEps{-}\wt\vv)
\Big)\,\d\xx\bigg)\d t+
\int_\varOmega\frac{|\pp_0|^2}{2\varrho_0}-\frac{|\pp_\EEps(T)|^2}{2\varrho_\EEps(T)}\,\d\xx
\\[-.5em]&=\nonumber
 \int_0^T\!\!\bigg(\bigg\langle\pdt{\pp_\EEps},\wt\vv-\vv_\EEps\bigg\rangle
 +\int_\varOmega\!\Big(\!\varrho_\EEps\GRAVITY_\EEps\Cdot(\vv_\EEps{-}\wt\vv)
-\big(\mathscr{T}_\lambda(\bm\xi,\Fe){-}\pp_\EEps{\otimes}\vv_\EEps\big)
\Colon\strain(\vv_\EEps{-}\wt\vv)
\\[-.7em]&\hspace{11.5em}
-\bbD\strain(\wt\vv)\Colon\strain(\vv_\EEps{-}\wt\vv)
-\HYPER|\nabla^2\wt\vv|^{p-2}\nabla^2\wt\vv\Vdots\nabla^2(\vv_\EEps{-}\wt\vv)
\Big)\,\d\xx\bigg)\d t\,,
\label{EUL-L-strong+++}\end{align}
where $\langle\cdot,\cdot\rangle$ denotes the duality between
$W^{2,p}(\varOmega;\R^3)^*$ and $W^{2,p}(\varOmega;\R^3)$ and where, for
the last equality, we used the calculus like \eq{EUL-L-strong-calc}
but for the continuous-in-time limit which holds as an equality. To
complete the Minty-trick arguments, we replace $\wt\vv$
by $\vv_\EEps\pm a\wt\vv$ with some $a>0$ and then pass $a\to0+$. Thus we
prove that $\vv_\EEps$ satisfies the momentum equation
\eq{Euler1-hypoplast-heter} in the weak sense \eq{def-ED-Ch5-momentum}.

The limit passage in \eq{EUL-L-visco-ED+3discr} is even simpler:
note that, the right-hand side of \eq{EUL-L-visco-ED+3discr}
converges only weakly when using the Green formula for
\begin{align}
&\nonumber
\int_\varOmega\big((\overlinevvtau\Cdot\nabla)\overlinexitau\big)\Colon\wt{\bm\xi}\,\d\xx=
-\int_\varOmega({\rm div}\overlinevvtau)\overlinexitau\Colon\wt{\bm\xi}
+\overlinexitau\Colon\big((\overlinevvtau\Cdot\nabla)\wt{\bm\xi}\big)\,\d\xx
\\[-.4em]&\nonumber\hspace{6em}
\to-\int_\varOmega({\rm div}\vv_\EEps)\bm\xi_\EEps\Colon\wt{\bm\xi}
+{\bm\xi}_\EEps\Colon\big((\vv_\EEps\Cdot\nabla)\wt{\bm\xi}\big)\,\d\xx
=\int_\varOmega\!
\big((\vv_\EEps\Cdot\nabla)\bm\xi_\EEps\big)\Colon\wt{\bm\xi}\,\d\xx\,
\end{align}
for any test function $\wt{\bm\xi}$ smooth. The same argument is to be used
also for the convective term in \eq{EUL-L-visco-ED+2discr}, while
the other terms in \eq{EUL-L-visco-ED+2discr} can be handled by 
\eq{Euler-small-converge-strong-E} with 
\eq{Euler-small-converge-strong+L}.

Altogether, we have proved that $(\varrho_\EEps,\vv_\EEps,\Fe_\EEps,\bm\xi_\EEps)$
is a weak solution to the truncated system
\begin{subequations}\label{ED-Euler-large-polyconvex-eps}\begin{align}
\label{ED-Euler-large0-polyconvex-eps}
&\pdt{\varrho_\EEps}=
-{\rm div}\,\pp_\EEps\ \ \text{ with }\ \pp_\EEps=\varrho_\EEps\vv_\EEps\,,\\
\nonumber
     &\pdt{\pp_\EEps}={\rm div}\big(\mathscr{T}_\lambda(\bm\xi_\EEps,\Fe_\EEps)
     +\DD_\EEps-\pp_\EEps{\otimes}\vv_\EEps\big)
     +\varrho_\EEps\GRAVITY
 \\[-.0em]
    &\hspace*{4em}\text{ with }\ \,\DD_\EEps=\bbD\strain(\vv_\EEps)
 -{\rm div}\,\mathfrak{H}_\EEps\,,
 \text{ where }\  
      \mathfrak{H}_\EEps=\HYPER|\nabla^2\vv_\EEps|^{p-2}\nabla^2\vv_\EEps\,,
\label{ED-Euler-large1-polyconvex-eps}
\\[-.3em]
&\pdt{\Fe_\EEps}=(\nabla\vv_\EEps)\Fe_\EEps
-\Fe_\EEps\mathscr{L}_\lambda(\bm\xi_\EEps,\Fe_\EEps)-(\vv_\EEps\Cdot\nabla)\Fe_\EEps\,,
\ \ \text{ and}
 \label{ED-Euler-large2-polyconvex-eps}
\\[-.1em]\label{ED-Euler-large3-polyconvex-eps}
&\pdt{\bm\xi_\EEps}=-(\vv_\EEps\Cdot\nabla)\bm\xi_\EEps\,,
\end{align}\end{subequations}
again with the boundary conditions \eq{Euler-small-BC-hyper} written
for $\TT=\mathscr{T}_\lambda(\bm\xi_\EEps,\Fe_\EEps)$.

Like \eq{ineq-for-sparsity}, we have now the sparsity equation
$\DT\sigma=({\rm div}\,\vv)\,\sigma$ for $\sigma=1/\varrho$ which
we now can test not only by $|\sigma|^{s-2}\sigma$ but also by
${\rm div}(|\nabla\sigma|^{r-2}\nabla\sigma)$. Analogously as
\eq{calulus-nonlin-rho-r}--\eq{est-of-rho-disc}, this test gives the estimate
\begin{align}\label{bound-for-sparsity}
\Big\|\frac1\varrho\Big\|_{L^\infty(I;W^{1,r}(\varOmega))}^{}\le C_r\ \ \text{ for any }\
1\le r<\infty\,.
\end{align}
In particular, for $r>3$ it yields that $\min\varrho>0$. This also improves
the estimate of $\vv$ in $L^\infty(I;L^a(\varOmega;\R^3))$ for $a>2$
which follows from \eq{EUL-L-est-v-2-eps+} as
\begin{align}
\|\vv\|_{L^\infty(I;L^2(\varOmega;\R^3))}^2
=\int_0^T\!\!\!\int_\varOmega\frac1\varrho\,\Big|\frac{\pp}{\sqrt\varrho}\Big|^2\,\d\xx\d t\le
\Big\|\frac1\varrho\Big\|_{L^\infty(I\times\varOmega)}^{}
\Big\|\frac{\pp}{\sqrt\varrho}\,\Big\|_{L^\infty(I;L^2(\varOmega;\R^3))}^{}\,,
\end{align}
which is bounded due to \eq{bound-for-sparsity} and the fist
estimate in \eq{Euler-small-est1}.

\medskip{\it Step 5: Elimination of the truncation}. Moreover, from
\eq{ED-Euler-large2-polyconvex-eps}, we also obtain an information about
$1/\det\Fe_\EEps$. Using the calculus 
\begin{align}\nonumber
\ \ \DT J&={\rm Cof}\Fe\Colon\DT\Fe=J\Fe^{-\top}\!\!\Colon\DT\Fe
=J\bbI\Colon\DT\Fe\Fe^{-1}\!=J\bbI\Colon\big(\Nabla\vv{-}\Fe\Lp\Fe^{-1}\big)
\\&=J{\rm div}\,\vv-J\Fe^{-\top}\!\Colon\Fe\Lp
=J{\rm div}\,\vv
-J\hspace{-1.8em}\lineunder{\bbI\Colon\Lp}{$\ \ \ \
={\rm tr}\Lp=0$}\hspace{-1.8em}=({\rm div}\,\vv)J\,,
\label{DT-det-extended}
\\[-2.3em]\nonumber\end{align}
we obtain the equation for $J_\EEps:=\det\Fe_\EEps$, namely
\begin{align}
\ \ \DT J_\EEps=({\rm div}\,\vv)J_\EEps\,.
\label{DT-det-extended+}\end{align}
From the assumption \eq{EUL-L-ass-IC} on the initial condition
$\det\Fezero$ and from the choice of $\lambda$ in Step~1, we have
$0<1/\det\Fe_\EEps<\lambda$ a.e. Moreover, recalling also the estimate
\eq{est-of-E-disc} for $r>3$ and again the choice of $\lambda$
in Step~1, we have also $|\Fe_\EEps|<\lambda$ a.e.\ on $I{\times}\varOmega$.
Thus we can omit the truncation $\varphi_\lambda^{}$ of $\varphi$ for the
obtained $\Fe$. Thus $(\varrho,\vv,\Fe,\bm\xi)$ is a solution to
\eq{Euler-hypoplast-p-heter} in the sense of Definition~\ref{def-ED-Ch5}.

\medskip{\it Step 6: Dissipation energy balance}.
It is now important that the tests and then all the subsequent calculations
leading to \eq{EUL-L-engr-balance} integrated over a current time
interval $[0,t]$ are really legitimate. The particular terms (including
time derivatives) are bounded in the norms indicated in \eq{EUL-L-est-dt}.

Testing \eq{Euler0-hypoplast-heter} by $|\vv|^2\in
L_{\rm w*}^{p/2}(I;L^\infty(\varOmega;\R^3))$, which leads to the kinetic energy,
can rely on that it is, for $p>3$, in
duality both with $\pdt{}\varrho\in L^p(I;L^r(\varOmega))$ and with
${\rm div}(\varrho\vv)=\vv\Cdot\nabla\varrho+\varrho({\rm div}\,\vv)
\in L^p(I;L^r(\varOmega))$ like \eq{EUL-L-est-p+}.

Furthermore, to execute the calculation \eq{Euler-hypoplast-test-momentum}
leading to the stored energy and the dissipation rate, we need to legitimate
the test of \eq{Euler1-hypoplast-heter} by $\vv\in
L^\infty(I;L^2(\varOmega;\R^3))\,\cap\,L^{p}(I;W^{2,p}(\varOmega;\R^3))$,
the test of $\DT{\Fe}\in L^p(I;L^r(\varOmega;\R^{3\times3}))$ in
\eq{Euler2-hypoplast-heter} by 
$\varphi_F'(\bm\xi,\Fe)\in L^\infty(I{\times}\varOmega;\R^{3\times3})$, and 
the test of $\DT{\bm\xi}\in L^p(I;L^r(\varOmega;\R^3))$ in
\eq{Euler3-hypoplast-heter} by
$\varphi_X'(\bm\xi,\Fe)\in L^\infty(I{\times}\varOmega;\R^3)$.
In the test of \eq{Euler1-hypoplast-heter}, we can rely on that
$(\vv\Cdot\nabla)\Fe\in L^p(I;L^r(\varOmega;\R^{3\times3})$ and
$(\nabla\vv)\Fe-\Fe\big[\zetap(\bm\xi,\cdot)^*\big]'\big({\rm dev}
(\Fee^{\!\!\top}\Se)\big)\in L^\infty(I{\times}\varOmega;\R^{3\times3})$.
\end{proof}

\begin{remark}[{\sl The assumptions (\ref{EUL-L-ass}a,b)}]
\upshape\label{rem-ass}
It should be pointed out that the truncation \eq{cutoff-stress} allows
also  for  wide
applicability of the qualification \eq{EUL-L-ass-zeta} of $\zetap$:
since $\zetap(X,\cdot)$ is minimized at the set $\{0\}$, 
the convex conjugate $\zetap(X,\cdot)^*$ is continuous and equals 0 at 0,
and also its derivative is zero at 0, which facilitates
the Lipschitz continuity in \eq{EUL-L-ass-zeta} when
$f:(X,F)\mapsto[\varphi_\lambda(X,\cdot)]'(F)$,
as used for \eq{EUL-L-heter-Gronwall+}  above.
\end{remark}

\vspace*{-.3em}

\begin{remark}[{\sl Higher regularity of $\bm\xi$}]\upshape
From \eq{ultimate-}, we have the equation
$\DT{\bm A}=-\bm A(\nabla\vv)$ for the so-called distortion
$\AA=\nabla\bm\xi$ with the initial condition $\AA|_{t=0}^{}=\bbI$.
Then we can obtain the estimate
$\|\AA\|_{L^\infty(I;W^{1,r}(\varOmega;\R^{3\times3}))}^{}\le C_r$, so that
$\|\bm\xi\|_{L^\infty(I;W^{2,r}(\varOmega;\R^{3}))}^{}\le C_r$ for any $r<\infty$.
\end{remark}

\begin{remark}[{\sl More general $\GRAVITY$}]\upshape
When having made a two-step convergence, we could allow for a weaker
qualification of $\GRAVITY$ than \eq{Euler-small-ass-g},
namely $\GRAVITY\in L^1(I;L^\infty(\varOmega;\R^3))$ and
perform the above time discretization with some approximation
of such $\GRAVITY$ in $L^\infty(I{\times}\varOmega;\R^3)$.
\end{remark}

\begin{remark}[{\sl A gradient on inelastic distortion rate}]\upshape
Sometimes, some gradient theories are used also in the
inelastic flow rule \eq{Euler3-hypoplast}. In the Eulerian model,
most straightforward variant is a dissipative gradient term as 
\begin{align}\label{Euler3-hypoplast+}
&\zetap(\bm\xi,\cdot)'(\Lp)-\ell\Delta\Lp={\rm dev}\MM\ \
\text{ with the Mandel stress }\
\MM=\Fee^{\!\!\top}\Se\,,
\end{align}
where $\ell>0$ is some length-scale parameter, to be completed by
the boundary condition $(\vv\Cdot\nabla)\Lp=\bm0$.
At small strain, this additional gradient leads to better propagation
of elastic waves, cf.\ \cite[Sect.\,3.4]{Roub24SGTL}. Here this
term augments the energy-dissipation balance \eq{EUL-L-engr-balance}
by the dissipation rate $\ell|\nabla\Lp|^2$ and allows for a certain weakening
of the Lipschitz-continuity assumption on $\zetap$ in \eq{EUL-L-ass-zeta}
because $\nabla\Lp$ can be estimated directly from the energy dissipation
balance and then use in \eq{EUL-L-heter-Gronwall+}.
On the other hand, it does not allow to eliminate
$\Lp$ by a pointwise manipulation as we did in \eq{Euler-hypoplast-p-heter}
and a separate time discretization of \eq{Euler3-hypoplast+}
should be done in the way that, instead of \eq{EUL-L-disc2},
we would consider a coupled system for $(\Fek,\Lp_{\tau}^k)$:
\begin{subequations}\begin{align}
&\!\!\frac{\Fek{-}\Fekk\!\!}\tau\,
=(\nabla\vv_\tau^k)\Fek-(\vv_\tau^k\Cdot\nabla)\Fek-\Fek\Lp_{\tau}^k\ \ \text{ and}
\label{EUL-L-disc2+}
\\[-.0em]
&\!\!\zetap(\bm\xi_{\tau}^k,\cdot)'
\big(\Lp_{\tau}^k\big)-\ell\Delta\Lp_{\tau}^k
={\rm dev}\MM_{\tau}^{k}\ \ \text{ with }\ \
\MM_{\tau}^{k}=(\Fek)^\top[\varphi_\lambda^{}]_F'(\bm\xi_{\tau}^{k},\Fek)\,.
\label{EUL-L-disc2++}
\end{align}\end{subequations}
\end{remark}

\begin{remark}[{\sl The symmetric hyperstress}]\upshape
One can consider a weaker hyperstress
$\mathfrak{H}=\HYPER|\nabla\strain(\vv)|^{p-2}\nabla\strain(\vv)$
instead of $\HYPER|\nabla^2\vv|^{p-2}\nabla^2\vv$. As in \cite{RouTom23IFST},
expressing the vorticity ${\rm curl}\,\strain(\vv)=\nabla\bm{w}$ with $\bm{w}$
being the axial vector of the skew-symmetric part of $\nabla\vv$, i.e.\
$\bm{w}{\times}\zz=(\frac12\nabla\vv-\frac12(\nabla\vv)^\top)\zz$ for
any $\zz$, one can see that the skew-symmetric part of $\nabla\vv$
and thus the whole $\nabla\vv$ is controlled by $\strain(\vv)$.
This weaker form of the hyperstress, making the Cauchy stress symmetric,
is physically even more consistent than $\mathfrak{H}$ in \eq{Euler1-hypoplast}.
\end{remark}

\begin{remark}[{\sl Discrete Gronwall inequality}]\upshape
For readers convenience, let us remind the following Gronwall-type
inequality 
\begin{equation}\label{9.2}
y_k\le\bigg(C+\tau\sum_{\ell=1}^{k-1}b_\ell\bigg)\,{\rm e}^{\tau\sum_{\ell=1}^{k-1}a_\ell}
\end{equation}
provided $y_k\le C+\tau\sum_{\ell=1}^{k-1}(a_\ell y_\ell+b_\ell)$ for any $k\ge0$; of
course, for $k=0$ it means $y_k\le C$. See e.g.\ 
\cite[Lemma~1.4.2]{QuaVal94NAPD}.
A more applicable variant needs the condition 
\begin{equation}\label{9.2+}
y_k\le C+\tau\sum_{\ell=1}^{k}\big(a_\ell y_\ell+b_\ell\big),
\end{equation}
for any $k\ge0$. Under the additional condition $\max_{1\le\ell\le k}a_\ell\le a$
for some $a$, it implies
$y_k\le (C+\tau\sum_{\ell=1}^{k-1}a_\ell y_\ell+\tau\sum_{\ell=1}^{k}b_\ell)/(1-\tau a)$
for any $k\ge0$. By \eq{9.2}, we then have 
\begin{equation}\label{9.2++}
y_k\le\frac{1}{1{-}a\tau}
\bigg(C+\tau\sum_{\ell=1}^kb_\ell\bigg)\,{\rm e}^{\tau\sum_{\ell=1}^ka_\ell/(1-a\tau)}
\ \ \ \mbox{ if }\ \tau<\frac1a\ \text{ with }\
a=\max_{1\le\ell\le k}a_\ell\,.
\end{equation}
Note that, in contrast to \eq{9.2}, the summation in \eq{9.2+} goes up to $k$,
which makes this condition more directly applicable than \eq{9.2}.
\end{remark}

\section{Extension towards some multiphysics models}\label{sec-multiphysics}

The possible arbitrary non-convexity of the stored energy, as allowed in
Sect.~\ref{sec-discretization}, allows for an extension that involves some
additional internal variables, such as damage (which may degrade especially
the stored energy) or some diffusant content (which may flow throughout the
viscoelastic medium via the Darcy or the Fick laws), or magnetization or
polarization in ferroic materials, etc. The truncated staggered
scheme can then be applied well. We will briefly illustrate this on two examples.

\subsection{Damage}

Let us consider the stored energy and the
viscoplastic potential depending on the scalar-valued parameter $\alpha$,
being interpreted as damage, whose evolution is governed
by the driving for{}ce  $\varphi_{\alpha}'$ through 
an additional dissipation, let us denote it by $\zetad$ depending
on $X\in\varOmega$ and on the rate of $\alpha$. Then, the system 
\eq{Euler-hypoplast} is to be enhanced as
\begin{subequations}\label{Euler-hypoplast-dam}
\begin{align}
\label{Euler0-hypoplast-dam}&\DT\varrho=-\varrho\,{\rm div}\,\vv\,,\\
\nonumber
     &\varrho\DT\vv={\rm div}(\TT{+}\DD)
+\varrho\GRAVITY\,\ \ \text{ with }\ 
\TT=\Se\Fee^{\!\!\top}\!\!+\varphi(\bm\xi,\Fe,\alpha)\bbI\,,\ \ \ \Se=\varphi_F'(\bm\xi,\Fe,\alpha)\,,
   \\[-.0em]
    &\hspace*{10.8em}\text{and }\ 
     \DD=\bbD\strain(\vv)-{\rm div}\,\mathfrak{H}\ \ \text{ with }\ 
\ \mathfrak{H}=\HYPER|\nabla^2\vv|^{p-2}\nabla^2\vv\,,
\label{Euler1-hypoplast-dam}
\\[-.4em]\label{Euler2-hypoplast-dam}
&\DT{\Fe}=
(\Nabla\vv)\Fe-\Fe\Lp\,,
\\\label{Euler3-hypoplast-dam}
&\zetap(\bm\xi,\alpha,\cdot)'(\Lp)={\rm dev}\MM\ \
\text{ with the Mandel stress }\
\MM=\Fee^{\!\!\top}\Se\,,
\\\label{Euler4-hypoplast-dam}
&
\partial\zetad(\bm\xi,\alpha,\cdot)(\DT\alpha)
\ni\varphi_{\alpha}'(\bm\xi,\Fe,\alpha)\,\
\text{ and}
\\[-.2em]\label{Euler5-hypoplast-dam}
&\DT{\bm\xi}=0\,.
\end{align}\end{subequations}
 Here $\zetad(\bm\xi,\alpha,\cdot)$ is admitted nonsmooth at $0$,
which turns \eq{Euler4-hypoplast-dam} into an inclusion in general and
then $\partial\zetad(\bm\xi,\alpha,\cdot)$ denotes the subdifferential. 
Naturally, the initial conditions \eq{IC-large} should be completed
by $\alpha|_{t=0}^{}=\alpha_0$.  The calculation
\eq{Euler-hypoplast-test-momentum} then expands 
by a term $\varphi_\alpha'(\bm\xi,\Fe,\alpha)\,\DT\alpha$. This term
can be substituted from \eq{Euler4-hypoplast-dam} tested by $\DT\alpha$,
which eventually expands  the {\it energy-dissipation balance}
\eq{EUL-L-engr-balance} by the damage dissipation rate as
\begin{align}\nonumber
\!\!\frac{\d}{\d t}\int_\varOmega\frac\varrho2|\vv|^2
+\varphi(\bm\xi,\Fe,\alpha)\,\d\xx
&+\!\int_\varOmega\Big(\bbD\strain(\vv)\Colon\strain(\vv){+}\HYPER|\nabla^2\vv|^p
\\[-.8em]&\qquad\
+\Lp\Colon\MM+\DT\alpha\,\zetad(\bm\xi,\alpha,\cdot)'(\DT\alpha)\Big)\,\d\xx
=\int_\varOmega\varrho\GRAVITY\Cdot\vv\,\d\xx\,.\ 
\label{EUL-L-energy-balance-dam}\end{align}

Assuming (as quite natural) convexity of $\zetad(X,\alpha,\cdot)$ with
minimum at $\DT\alpha=0$, we can rewrite the damage flow rule
\eq{Euler4-hypoplast-dam} in the form 
\begin{align}\label{damage-transformed}
\DT\alpha=\big[\zetad(\bm\xi,\alpha,\cdot)^*\big]'\big(\varphi_{\alpha}'(\bm\xi,\Fe,\alpha)\big)\,.
\end{align}
 Note that the convex conjugate $\zetad(\bm\xi,\alpha,\cdot)^*$ is
(assumed) continuously differentiable even if $\zetad(\bm\xi,\alpha,\cdot)$
may be not smooth. Thus we admit (although not required) a unidirectional
damage.  This form can be advantageously used for the time discretization.
The truncation \eq{cutoff-stress} of $\varphi$ is to be straightforwardly
modified by adding the variable $\alpha$ and then, like \eq{truncated},
we define
\begin{subequations}\label{truncated+}\begin{align}\label{stress-truncated++}
&\mathscr{T}_\lambda^{}(X,F,\alpha):=
[\varphi_\lambda]_F'(X,F,\alpha)F^\top\!{+}\varphi_\lambda^{}(X,F,\alpha)\,\bbI\,,
\\&\label{Lp-truncated-heter+}
\mathscr{L}_\lambda^{}(X,F,\alpha):=\big[\zetap(X,\alpha,\cdot)^*\big]'\big({\rm dev}(F^\top\![\varphi_\lambda^{}]_F'(X,F,\alpha))\big)\,,\ \ \ \text{ and }\  
\\&\label{damage-truncated-heter+}
\mathscr{D}_\lambda^{}(X,F,\alpha):=
\big[\zetad(X,\alpha,\cdot)^*\big]'\big([\varphi_\lambda^{}]_\alpha'(X,F,\alpha)\big)\,.
\end{align}\end{subequations}
The staggered scheme \eq{EUL-L-disc} can then be enhanced as 
\begin{subequations}\label{EUL-L-viscoelastodyn+disc-dam}
\begin{align}\label{EUL-L-viscoelastodyn+0disc-dam}
&\!\!\frac{\varrho_\etau^k{-}\varrho_\etau^{k-1}\!\!}\tau\,=
-{\rm div}\,\pp_\etau^k\ \ \ \text{ with }\ \
\pp_\etau^k=\varrho_\etau^k\vv_\etau^k\,,\!
\\[-.2em]&\nonumber
\!\!\frac{\pp_\etau^k{-}\pp_\etau^{k-1}\!\!}\tau\,=
{\rm div}\Big(
\mathscr{T}_\lambda(\bm\xi_{\etau}^{k-1},\Fekk,\alpha_{\eetau}^{k-1})
+\DD_\etau^k-\pp_\etau^k{\otimes}\vv_\etau^k\Big)+\varrho_\etau^k\GRAVITY_{\tau}^k\,,
\\[-.2em]&
\hspace*{4.5em}\text{ where }\ \ \DD_\etau^k=\bbD\strain(\vv_\etau^k)
-{\rm div}\,\mathfrak{H}_\etau^k\ \ \text{ with }\ \
\mathfrak{H}_\etau^k=\HYPER\big|\nabla^2\vv_\etau^k\big|^{p-2}\nabla^2\vv_\etau^k\,,\!
\label{EUL-L-viscoelastodyn+1disc-dam}
\\[-.1em]
&\!\!\frac{\Fek{-}\Fekk\!\!}\tau\,
=(\nabla\vv_\etau^k)\Fek-
\Fek\mathscr{L}_\lambda^{}(\bm\xi_{\etau}^{k},\Fek,\alpha_{\eetau}^k)
-(\vv_\etau^k\Cdot\nabla)\Fek\,,
\label{EUL-L-viscoelastodyn+2disc-dam}
\\[-.1em]
\label{EUL-L-viscoelastodyn+3disc-dam}
&\!\!\frac{\alpha_{\eetau}^k{-}\alpha_{\eetau}^{k-1}\!\!}\tau\,
=\mathscr{D}_\lambda^{}(\bm\xi_{\etau}^{k},\Fek,\alpha_{\eetau}^k)
-\vv_\etau^k\Cdot\nabla\alpha_{\eetau}^k\,,\ \ \text{ and}
\\[-.1em]
&\!\!\frac{\bm\xi_{\etau}^k{-}\bm\xi_{\etau}^{k-1}\!\!}\tau\,
=-(\vv_\etau^k\Cdot\nabla)\bm\xi_{\etau}^k\,.
\label{EUL-L-viscoelastodyn+4disc-dam}
\end{align}\end{subequations}
Note that, like \eq{EUL-L-disc}, also this system is decoupled
(staggered) in the sense that (\ref{EUL-L-viscoelastodyn+disc-dam}a,b)
gives $(\varrho_\etau^k,\pp_\etau^k)$ and thus $\vv_\etau^k$ which is then used
for \eq{EUL-L-viscoelastodyn+4disc-dam} to obtain $\bm\xi_{\etau}^k$
and then for (\ref{EUL-L-viscoelastodyn+disc-dam}c,d)
to obtain $(\Fek,\alpha_\etau^k)$.

The analysis from Sect.\,\ref{sec-discretization} can be extended
quite straightforwardly, starting again from the energy imbalance
\eq{EUL-L-basic-engr-balance-disc} now naturally with 
$\mathscr{T}_\lambda^{}(\bm\xi_{\eetau}^{m-1},\Femm,\alpha_{\eetau}^{m-1})$
in place of $\mathscr{T}_\lambda^{}(\bm\xi_{\eetau}^{m-1},\Femm)$. In the limit
passage, we now need also the strong convergence of $\overlinealphaetau$
 and of $\underlinealphaetau$, i.e.\ of the piecewise-constant
interpolants of $\{\alpha_{\eetau}^k\}_{k=0}^{T/\tau}$. To this goal, we aim to
get an estimate of $\nabla\overlinealphaetau$. It is important  that 
\eq{EUL-L-viscoelastodyn+3disc-dam} is explicit in terms
of the time difference thanks to the transformation \eq{damage-transformed}
so that it can be tested by
${\rm div}(|\nabla\alpha_{\eetau}^k|^{r-2}\nabla\alpha_{\eetau}^k)$.
This test needs the assumption 
\begin{align}\label{damage-transformed-ass}
\forall \lambda>1:\ \ \ (X,F,\alpha)\mapsto
\big[\zetad(X,\alpha,\cdot)^*\big]'\big([\varphi_\lambda^{}]_{\alpha}'(X,F,\alpha)\big)
\ \text{ is Lipschitz continuous}\,.
\end{align}
 Like \eq{calulus-nonlin-rho-r}--\eq{calulus-nonlin-rho-r+}, 
this then gives the inequality
\begin{align}\nonumber
&
\int_\varOmega\!\frac{|\nabla\alpha_{\eetau}^k|^{ r}\!-|\nabla\alpha_{\eetau}^{k-1}|^{ r}\!\!}{{ r}\tau}\,\d\xx
\le\!\int_\varOmega\!|\nabla\alpha_{\eetau}^k|^{r-2}_{^{^{}}}\Big(
\nabla\mathscr{D}_\lambda^{}(\bm\xi_{\etau}^{k},\Fek,\alpha_{\eetau}^k)\Cdot
\nabla\alpha_{\eetau}^k\!
-\nabla\big((\vv_{\eetau}^k\Cdot\nabla)\alpha_{\eetau}^k
\big)\Colon\nabla\alpha_{\eetau}^k\Big)\,\d\xx
\\[-.4em]&\nonumber\ 
=\int_\varOmega|\nabla\alpha_{\eetau}^k|^{r-2}_{^{^{}}}\bigg(
\big[\mathscr{D}_\lambda^{}\big]_X'(\bm\xi_{\etau}^{k},\Fek,\alpha_{\eetau}^k)
\nabla\bm\xi_{\etau}^{k}\Cdot\nabla\alpha_{\eetau}^k
+\big[\mathscr{D}_\lambda^{}\big]_{F}'(\bm\xi_{\etau}^{k},\Fek,\alpha_{\eetau}^k)\nabla\Fek\Cdot\nabla\alpha_{\eetau}^k
\\[-.6em]&\nonumber\qquad\qquad
+\big[\mathscr{D}_\lambda^{}\big]_{\alpha}'(\bm\xi_{\etau}^{k},\Fek,\alpha_{\eetau}^k)\nabla\alpha_{\eetau}^k\Cdot\nabla\alpha_{\eetau}^k\!
+\frac{\!{\rm div}\,\vv_{\etau}^{k}\!}{ r}\,|\nabla\alpha_{\eetau}^k|^2\!
-
(\nabla\alpha_{\eetau}^k{\otimes}\nabla\alpha_{\eetau}^k)\Colon\strain(\vv_{\eetau}^k)\bigg)\,\d\xx
\\[-.4em]&
\ \le\ell_\lambda^{ r}\!
+\|\nabla\bm\xi_{\etau}^{k}\|_{L^{ r}(\varOmega;\R^{3\times3})}^{ r}\!
+\|\nabla\Fek\|_{L^{ r}(\varOmega;\R^{3\times3\times3})}^{ r}
+\!\Frac32\|\strain(\vv_{\eetau}^k)\|_{L^\infty(\varOmega;\R^{3\times3})}^{ r}
+\|\nabla\alpha_{\eetau}^k\|_{L^{ r}(\varOmega;\R^3)}^{ r}
\label{est-of-nabla-alpha}\end{align}
with $\ell_\lambda^{}$ denoting the Lipschitz constant from the condition
\eq{damage-transformed-ass}. Here we used the calculus based on the Green
formula as \eq{calulus-nonlin-rho-r} now for 
the scalar variable $\alpha$  instead of $\varrho$. Furthermore, 
the $\alpha$-dependence of \eq{EUL-L-viscoelastodyn+2disc-dam}
gives rise to an additional term
$\int_\varOmega(\Fek[\mathscr{L}_\lambda^{}]_\alpha'(\bm\xi_{\etau}^{k},\Fek,\alpha_{\etau}^{k}))
\nabla\alpha_{\etau}^{k})$ $\Vdots(|\nabla\Fek|^{r-2}\nabla\Fek)\,\d\xx$
in \eq{EUL-L-heter-Gronwall+}. This term can be estimated like the
$\bm\xi$-term in \eq{EUL-L-heter-Gronwall+}, which gives the terms
as in \eq{est-of-nabla-alpha}. Therefore, \eq{EUL-L-heter-Gronwall+} with
this additional term must be estimated jointly with \eq{est-of-nabla-alpha}.

Using the already available estimates \eq{Euler-small-est1},
\eq{est-of-xi-disc}, and \eq{est-of-E-disc} and the discrete Gronwall
inequality, for  any  sufficiently small time steps $\tau>0$, we obtain the
estimate
\begin{align}\label{est-damage-grad}
\|\nabla\overlinealphaetau\|_{L^\infty(I;L^{ r}(\varOmega;\R^3))}^{}\le C
\end{align}
 together with \eq{est-of-E-disc}; specifically,
\eq{est-of-nabla-alpha} suggests $\tau>0$ small enough
as used for \eq{est-of-xi-disc} and \eq{est-of-E-disc}. 
In particular, 
this implies the strong convergence $\overlinealphaetau\to\alpha$
 and also $\underlinealphaetau\to\alpha$ 
in $L^2(I{\times}\varOmega)$ due to a (generalized) Aubin-Lions theorem when
taking into account also the information about
$\pdt{}\alpha_\tau\in L^2(I{\times}\varOmega)$ available by comparison
from \eq{EUL-L-viscoelastodyn+3disc-dam}.

Let us note that the Lipschitz continuity assumption 
\eq{damage-transformed-ass} needs, in particular,
$[\zetad(X,\alpha,\cdot)^*]'$ Lipschitz continuous, which
requires that $\zetad(X,\alpha,\cdot)$ has at least quadratic growth
uniformly in $(X,\alpha)$. This excludes the rate-independent damage but,
on the other hand, still admits a unidirection{}al  damage (i.e.\ with
healing  forbidden  so that always $\DT\alpha\le0$) as typically used
in engineering models. Damage models typically assume
$0\le\alpha\le1$. In the damage without healing, this can be
 ensured  by $\varphi_\alpha'(X,F,0)\ge0$ if the initial condition
satisfies $0\le\alpha_0\le1$. If healing is allowed (as often  e.g.\ 
in  some  geophysical  or biological  models), then
also $\varphi_\alpha'(X,F,1)=0$ is needed.

The Definition~\ref{def-ED-Ch5} of the weak solution is straightforwardly
modified by considering the $\alpha$-dependence with
$\alpha\in L^\infty(I;W^{1,r}(\varOmega))\,\cap\,H^1(I;L^2(\varOmega))$ valued in
$[0,1]$ and \eq{damage-transformed} holding a.e.\ on $I{\times}\varOmega$
together with the initial condition $\alpha|_{t=0}^{}=\alpha_0$. The above sketched
arguments augmenting the proof of Proposition~\ref{prop-ED-Ch5-existence} leads to:

\begin{proposition}
Let (\ref{EUL-L-ass}b-e)
with $\zetap$ modified as $\alpha$-dependent hold for $p>r>3$ with
\begin{subequations}\label{EUL-L-ass-dam}
\begin{align}\nonumber
&\!\varphi\in C^1(\barOmega{\times}\GL_3^+{\times}[0,1]), \ \
\forall (X,F){\in}\varOmega{\times}\GL_3^+:\ \ \varphi(X,F,\cdot)\
\text{nondecreasing and}
\\&
\hspace{18em}
\inf_{X\in\varOmega,\ F\in \GL_3^+,\ \alpha\in[0,1]}\varphi(X,F,\alpha)>-\infty\,,
\label{EUL-L-ass-dam++}
\\&\nonumber
\zetad\in \varOmega{\times}[0,1]{\times}\R\to[0,+\infty]\ \ \text{ satisfies }\
\eq{damage-transformed-ass}\,, \text{ and }
\\&\nonumber
\qquad\qquad
\exists\delta>0\ \forall(X,\alpha,\DT\alpha)\in\varOmega{\times}[0,1]{\times}\R:\ \ 
\zetad(X,\alpha,\DT\alpha)\ge\tfrac12\delta|\DT\alpha|^2
\\[-.1em]&\hspace{4em}\nonumber
\forall (X,\alpha){\in}\varOmega{\times}[0,1]:\ \ \
\zetad(X,\alpha,\cdot)\ \text{is convex with a minimum at }\{0\}\,, \text{ and }
\\&\hspace{13.5em}
\zetad(X,\alpha,\DT\alpha)|_{\DT\alpha>0}^{}\!=+\infty
\ \text{ or }\ \varphi_\alpha'(X,F,1)=0\,,
\label{EUL-L-ass-zeta+}
\end{align}\end{subequations}
Then the discrete scheme \eq{EUL-L-viscoelastodyn+disc-dam}
gives approximate solutions which converge in terms of subsequences to 
weak solutions to the problem \eq{Euler-hypoplast-dam} with
the boundary conditions \eq{Euler-small-BC-hyper} and
the initial conditions \eq{IC-large} together with 
the condition $\alpha|_{t=0}^{}=\alpha_0\in W^{1,r}(\varOmega)$ with
$0\le\alpha_0\le1$ on $\varOmega$.
\end{proposition}

Let us also emphasize that, due to the hyperviscosity, we did not need
to include any gradient of damage in the model
 to obtain \eq{est-damage-grad} and also we can admit a complete damage
in the conservative part, i.e.\ $\varphi(X,F,0)\equiv0$ is possible. Anyhow,
involving a term ${\rm div}(\ell^r|\nabla\alpha|^{r-2}\nabla\alpha)$
in the right-hand side of \eq{Euler4-hypoplast-dam} with some length scale
parameter $\ell>0$ is  possible and  often a desirable modeling aspect.
In convective models, this 
leads to a (symmetric) Korteweg-type contribution to the Cauchy stress, namely
$\ell^r(|\nabla\alpha|^{r-2}\nabla\alpha{\otimes}\nabla\alpha-\frac1r|\nabla\alpha|^r\bbI)$.

The weak formulation of the inclusion \eq{Euler4-hypoplast-dam} with the
augmented right-hand side by ${\rm div}(\ell^r|\nabla\alpha|^{r-2}\nabla\alpha)$
leads to the variational inequality
\begin{align}
&\nonumber
\int_0^T\!\!\!\int_\varOmega\zetad(\bm\xi,\alpha,\DT\alpha)\,\d\xx\d t
\le\int_\varOmega\frac{\ell^r}r|\nabla\alpha(T)|^r
-\frac{\ell^r}r|\nabla\alpha_0|^r\,\d\xx
+\int_0^T\!\!\!\int_\varOmega\bigg(\zetad(\bm\xi,\alpha,\wt\alpha)
\\[-.4em]&\quad\ +
\varphi_{\alpha}'(\bm\xi,\Fe,\alpha)
(\wt\alpha{-}\DT\alpha)
-\ell^r|\nabla\alpha|^{r-2}\nabla\alpha\Cdot\nabla\wt\alpha
-|\nabla\alpha|^{r-2}(\nabla\alpha{\otimes}\nabla\alpha-\frac1r|\nabla\alpha|^r\bbI)\Colon\strain(\vv)\bigg)\,\d\xx\d t
\nonumber\end{align}
for any $\wt\alpha\in L^\infty(I;W^{1,r}(\varOmega))$; here
we substituted $\int_\varOmega|\nabla\alpha|^{r-2}\nabla\alpha\Cdot\nabla\DT\alpha\,\d\xx
=\int_\varOmega|\nabla\alpha|^{r-2}\nabla\alpha\Cdot\nabla(\pdt{}\alpha+\vv\Cdot\nabla\alpha)\,\d\xx
=\int_\varOmega\frac1r\pdt{}|\nabla\alpha|^r
+|\nabla\alpha|^{r-2}(\nabla\alpha{\otimes}\nabla\alpha-\frac1r|\nabla\alpha|^r\bbI)\Colon\strain(\vv)\,\d\xx$ while using the
calculus \eq{calulus-nonlin-rho-r} now again for $\alpha$ instead of $\varrho$. 
The staggered scheme \eq{EUL-L-viscoelastodyn+disc-dam} is then to
be modified. Specifically, (\ref{EUL-L-viscoelastodyn+disc-dam}a,c,d)
are to be completed by the equations
(\ref{EUL-L-viscoelastodyn+disc-dam}b,d) enhanced as
\begin{subequations}\label{EUL-L-viscoelastodyn+disc-dam+}
\begin{align}
&\nonumber
\!\!\frac{\pp_\etau^k{-}\pp_\etau^{k-1}\!\!}\tau\,=
{\rm div}\Big(
\mathscr{T}_\lambda(\bm\xi_{\etau}^{k-1},\Fekk,\alpha_{\eetau}^{k-1})
+\DD_\etau^k-\pp_\etau^k{\otimes}\vv_\etau^k
\\[-.6em]&
\hspace*{10.5em}
+\ell^r\big(|\nabla\alpha|^{r-2}\nabla\alpha_{\etau}^{k}{\otimes}\nabla\alpha_{\etau}^{k}-\frac1r|\nabla\alpha_{\etau}^{k}|^r\bbI\big)\Big)
+\varrho_\etau^k\GRAVITY_{\tau}^k\ \ \text{ and}
\label{EUL-L-viscoelastodyn+1disc-dam+}
\\[-.3em]\nonumber
&\!\!
\partial\zetad(\bm\xi_{\eetau}^k,\alpha_{\eetau}^k,\cdot)\Big(
\frac{\alpha_{\eetau}^k{-}\alpha_{\eetau}^{k-1}\!\!}\tau
+\vv_\etau^k\Cdot\nabla\alpha_{\eetau}^k\Big)
\\[-.3em]&
\hspace*{11.5em}
\ni[\varphi_\lambda]_{\alpha}'(\bm\xi_{\eetau}^{k-1},\Fekk,\alpha_{\eetau}^{k-1})
+{\rm div}\big(\ell^r|\nabla\alpha_{\etau}^{k}|^{r-2}\nabla\alpha_{\etau}^{k}\big)
\label{EUL-L-viscoelastodyn+3disc-dam+}
\end{align}\end{subequations}
with $\DD_\etau^k$ from \eq{EUL-L-viscoelastodyn+1disc-dam+}. Of course,
the boundary conditions \eq{BC-disc} should use the correspondingly
augmented Cauchy stress and, moreover, 
\eq{EUL-L-viscoelastodyn+3disc-dam+} needs an additional
boundary condition, mostly considered as $\nn\Cdot\nabla\alpha_{\etau}^{k}=0$
on $\varGamma$. Notably, these additional gradient terms couple
the equations for $(\varrho_{\eetau}^k,\pp_{\eetau}^k,\Fek,\alpha_{\eetau}^k)$
while the equation for $\bm\xi_{\eetau}^k$ still remains decoupled.

The estimate \eq{est-damage-grad} can be then obtained directly
from the energy dissipation balance which is then enhanced by the
contribution $\frac1r\ell^r|\nabla\alpha|^r$ to the stored energy.
More in detail, testing \eq{EUL-L-viscoelastodyn+3disc-dam+} by
the discrete convective derivative
$\frac{\alpha_{\eetau}^k{-}\alpha_{\eetau}^{k-1}\!\!}\tau+\vv_\etau^k\Cdot\nabla\alpha_{\eetau}^k$ uses the calculus as in \eq{calulus-nonlin-rho-r} once again for $\alpha_{\etau}^{k}$
instead of $\varrho_{\etau}^{k}$, namely 
\begin{align}\nonumber
\int_\varOmega
{\rm div}\big(\ell^r|\nabla\alpha_{\etau}^{k}|^{r-2}\nabla\alpha_{\etau}^{k}\big)
\big(\vv_\etau^k\Cdot\nabla\alpha_{\eetau}^k\big)\,\d\xx
&=\int_\varOmega
\ell^r|\nabla\alpha_{\etau}^{k}|^{r-2}\Big(\frac1r|\nabla\alpha_{\etau}^{k}|^2\bbI-
\nabla\alpha_{\etau}^{k}{\otimes}\nabla\alpha_{\etau}^{k}\Big)\Colon\strain(\vv_\etau^k)\,\d\xx\,,
\end{align}
which then cancels with the Korteweg-like stress in
\eq{EUL-L-viscoelastodyn+1disc-dam+}.
Therefore, instead of \eq{EUL-L-basic-engr-balance-disc}, we now
obtain 
\begin{align}\nonumber
\!\!\!\!\!&\int_\varOmega\!\frac{|\pp_\etau^k|^2\!\!}{2\varrho_\etau^k\!\!\!}
+\frac{\ell^r}r|\nabla\alpha_\etau^k|^r\d\xx
+\tau\!\sum_{m=1}^k\int_\varOmega\!
\bbD\strain(\vv_\etau^m)\Colon\strain(\vv_\etau^m)
+\HYPER|\nabla^2\vv_\etau^m|^p
+\delta\Big|\frac{\alpha_{\eetau}^m{-}\alpha_{\eetau}^{m-1}\!\!\!\!}\tau+\vv_\etau^m\Cdot\nabla\alpha_{\eetau}^m\Big|^2\d\xx
\\[-.5em]&\nonumber\hspace{2em}
\le\int_\varOmega\!\frac{|\pp_0|^2\!\!}{2\varrho_0\!\!}
+\frac{\ell^r}r|\nabla\alpha_0|^r\,\d\xx
+\tau\!\sum_{m=1}^k\int_\varOmega\!\varrho_\etau^m\GRAVITY_{\tau}^m\Cdot\vv_\etau^m\!
-\mathscr{T}_\lambda^{}(\bm\xi_{\eetau}^{m-1},\Femm,\alpha_{\eetau}^{m-1})\Colon\strain(\vv_\etau^m)
\\[-.5em]&\hspace{12em}
+[\varphi_\lambda]_{\alpha}'(\bm\xi_{\eetau}^{m-1},\Femm,\alpha_{\eetau}^{m-1})
\Big(\frac{\alpha_{\eetau}^m{-}\alpha_{\eetau}^{m-1}\!\!\!\!}\tau+\vv_\etau^m\Cdot\nabla\alpha_{\eetau}^m\Big)
\,\d\xx\,,
\label{EUL-L-basic-engr-balance-disc-dam}\end{align}
where $\delta>0$ is from \eq{EUL-L-ass-zeta+}.
The last term can be absorbed in the left-hand side, which gives additionally
the estimate of $\pdt{}\alpha_\tau+\overlinevvtau\Cdot\nabla\overlinealphaetau$
in $L^2(I{\times}\varOmega)$ and, using also \eq{EUL-L-est-v-2-eps++} and
\eq{est-damage-grad}, then also of $\pdt{}\alpha_\tau$ in
$L^2(I{\times}\varOmega)$. Due to the nonlinear Korteweg stress, we now need 
the strong convergence of $\nabla\overlinealphaetau\to\nabla\alpha$.
This  can be obtained by the strong monotonicity of this
extra term $-{\rm div}(\ell^r|\nabla\alpha|^{r-2}\nabla\alpha)$
in \eq{Euler4-hypoplast-dam}, cf.\ 
\cite{RouTom21CMPE} for a linearized convected model discretized
by a fully coupled (non-staggered) scheme.

\begin{remark}[{\sl Phase-field fracture}]\upshape
Often, damage concentrates on the small regions, i.e.\ surfaces
representing cracks. In the context of the damage concept, this can be modeled by
a ``penalization'' term in the stored energy $\varphi$ forcing the 3-dimensional
damaged regions to have a big energy. Together with the gradient term as in
\eq{EUL-L-viscoelastodyn+3disc-dam+}, this approximative technique leads to
the so-called phase-field fracture. Damage or fracture processes can be very
fast (sometimes even modeled as rate independent which allows very big
or even unlimited speed of fracture), which may lead to emission of waves
and then the combination with inertia is particularly important.
\end{remark}

\subsection{Diffusion by the chemical potential gradient}\label{sec-diffuse}

Another interesting internal variable, again denoted by $\alpha$
and expanding the stored energy $\varphi=\varphi(X,F,\alpha)$, is a
{\it diffusant concentration} considered as an {\it intensive variable},
valued in the interval $[0,1]$. Its evolution is governed by the
diffusive flux through a mobility coefficient $M=M(X,\alpha)$
multiplying the gradient of the chemical potential
$\mu=\varphi_\alpha'(\bm\xi,\Fe,\alpha)$. This may cover both the Darcy and
the Fick laws. Instead of \eq{Euler4-hypoplast-dam}, we thus consider the
balance equation
\begin{align}\label{EUL-L-diff}
\DT\alpha={\rm div}\big(M(\bm\xi,\alpha)\nabla\mu\big)
\ \ \ \text{ with }\ \mu\in\varphi_\alpha'(\bm\xi,\Fe,\alpha)+N_{[0,1]}(\alpha)\,.
\end{align}
Here $N_{[0,1]}(\cdot)$ denotes the set-valued normal-cone mapping to
the interval $[0,1]$, i.e.\
$$
N_{[0,1]}(\alpha)=\begin{cases}\quad\ 0&\text{for }\ 0<\alpha<1\,,\\[-.2em]
(-\infty,0]&\text{for }\ \alpha=0\,,\\[-.2em]
[0,+\infty)&\text{for }\ \alpha=1\,,\\[-.2em]
\quad\ \emptyset&\text{for }\ \alpha\not\in[0,1]\,;\end{cases}
$$
notably, $N_{[0,1]}(\cdot)=\partial\delta_{[0,1]}(\cdot)$ is the subdifferential
of the convex indicator function $\delta_{[0,1]}(\cdot)$ to the interval $[0,1]$.
This set-valued mapping $N_{[0,1]}$ thus causes \eq{EUL-L-diff} as an
inclusion and ensures that $\alpha$ is valued in $[0,1]$.

The calculation \eq{Euler-hypoplast-test-momentum} is again expanded 
by the term $\varphi_\alpha'(\bm\xi,\Fe,\alpha)\,\DT\alpha$, which
is now handled by testing \eq{EUL-L-diff} by $\mu$.
 Prescribing still the boundary condition $\nn\Cdot\nabla\mu=0$ on $\varGamma$,
 this gives the {\it energy-dissipation balance} 
\eq{EUL-L-engr-balance} enhanced by the diffusion
dissipation rate as
\begin{align}\nonumber
\!\!\frac{\d}{\d t}\int_\varOmega&\,\frac\varrho2|\vv|^2
+\varphi(\bm\xi,\Fe,\alpha)\,\d\xx
\\[-.5em]&\
+\!\int_\varOmega\bbD\strain(\vv)\Colon\strain(\vv)+\HYPER|\nabla^2\vv|^p
+\Lp\Colon\MM+M(\bm\xi,\alpha)|\nabla\mu|^2\,\d\xx
=\!\int_\varOmega\!\varrho\,\GRAVITY\Cdot\vv\,\d\xx\,.\ 
\label{EUL-L-energy-balance-diff}\end{align}
From \eq{EUL-L-energy-balance-diff}, one can see the estimate of
$\nabla\mu$ in $L^2(I{\times}\varOmega;\R^3)$. In contrast to the previous
damage model, the flow rule \eq{EUL-L-diff} thus allows directly
for estimation of $\nabla\alpha$ if
$\varphi$ is twice differentiable with $\varphi(X,F,\cdot)$ is strongly convex.
This can be seen from the calculus 
\begin{align}\nonumber
\nabla\mu&=\nabla
\big(\varphi_\alpha'(\bm\xi,\Fe,\alpha)+\partial\delta_{[0,1]}^{}(\alpha)\big)
\\&=\varphi_{X\alpha}''(\bm\xi,\Fe,\alpha)\Cdot\nabla\bm\xi
+\varphi_{F\alpha}''(\bm\xi,\Fe,\alpha)\Colon\nabla\Fe
+\varphi_{\alpha\alpha}''(\bm\xi,\Fe,\alpha)\nabla\alpha
+\partial^2\delta_{[0,1]}(\alpha)\nabla\alpha\,.
\label{nabla-mu}\end{align}
Using the positive semi-definiteness of  the generalized Hessian 
$\partial^2\delta_{[0,1]}(\alpha)$  of the convex indicator
function $\delta_{[0,1]}(\cdot)$ 
so that $\partial^2\delta_{[0,1]}(\alpha)(\nabla\alpha,\nabla\alpha)\ge0$,
we can test \eq{nabla-mu} by $\nabla\alpha$ and then obtain the inequality
\begin{align}\label{est-of-nabla-alpha+}
|\nabla\alpha|^2\le
\Big(\frac{\nabla\mu}{\varphi_{\alpha\alpha}''(\bm\xi,\Fe,\alpha)}
-\frac{\varphi_{X\alpha}''(\bm\xi,\Fe,\alpha)}{\varphi_{\alpha\alpha}''(\bm\xi,\Fe,\alpha)}
\Cdot\nabla\bm\xi
-\frac{\varphi_{F\alpha}''(\bm\xi,\Fe,\alpha)}{\varphi_{\alpha\alpha}''(\bm\xi,\Fe,\alpha)}
\Colon\nabla\Fe
\Big)\Cdot\nabla\alpha\,;
\end{align}
actually the handling of the generalized  Hessian  $\partial^2$ of the
convex nonsmooth function is rather formal (though standard) and
can be legitimized by mollifying $\delta_{[0,1]}(\cdot)$.
Due to the uniform positive definiteness of $\varphi(X,F,\cdot)$ ensuring
$\inf\varphi_{\alpha\alpha}''>0$ and due to the already available estimates
on $\nabla\mu$, $\nabla\bm\xi$, and $\nabla\Fe$, one can infer the
bound for $\nabla\alpha$ in $L^2(I{\times}\varOmega;\R^3)$. Together with an
information about $\pdt{}\alpha\in L^2(I;W^{1,2}(\varOmega)^*)$ obtained by
comparison from the first relation in \eq{EUL-L-diff}, it gives
compactness of $\alpha$'s in $L^2(I{\times}\varOmega)$, which gives a base
for the analysis of the time-discrete scheme which now should use, instead of
\eq{EUL-L-viscoelastodyn+3disc-dam}, 
\begin{align}\nonumber
\!\frac{\!\alpha_{\tau}^k{-}\alpha_{\tau}^{k-1}\!\!\!}\tau
={\rm div}&\big(M(\bm\xi_{\tau}^k,\alpha_{\tau}^k)\nabla\mu_{\tau}^k\big)
-\vv_\tau^k\Cdot\nabla\alpha_{\tau}^k
\\[-.4em]&\qquad
\ \,\text{ with }\
\mu_{\tau}^k\!\in[\varphi_\lambda]_\alpha'(\bm\xi_{\tau}^k,\Fek,\alpha_{\tau}^k)
+N_{[0,1]}(\alpha_{\tau}^k).\!
\label{EUL-L-viscoelastodyn+3disc-dif}
\end{align}

Here, for estimates, we need to assume
\begin{subequations}\label{ass-difuse}
\begin{align}\label{ass-difuse1}
&\zetap=\zetap(X,L)\ \ \ \text{ and }\ \ {\rm dev}(F^\top\varphi(X,\cdot,\alpha)'(F))
\text{ independent of }\alpha,
\text{ and}
\\&\sup_{X\in\varOmega,\ F\in\R^{3\times3},\ \alpha\in[0,1]}\!\!\!
\frac{|[\varphi_\lambda]_{(X,F),\alpha}''(X,F,\alpha)|}
{[\varphi_\lambda]_{\alpha\alpha}''(X,F,\alpha)}
+\frac1{\sqrt{M(X,\alpha)}[\varphi_\lambda]_{\alpha\alpha}''(X,F,\alpha)}
<\infty\,
\label{ass-difuse2}\end{align}
\end{subequations}
for any $\lambda>1$ fixed. 
Notably, \eq{ass-difuse1} says that the concentration $\alpha$ of the
diffusant does not influence the isochoric part of the model.
This is natural in isotropic media where swelling or squeezing act on
the volumetric part only. In particular, it makes 
$\mathscr{L}_\lambda^{}$ independent of $\alpha$.

From the equations (\ref{EUL-L-viscoelastodyn+disc-dam}a,b), we again obtain
\eq{EUL-L-basic-engr-balance-disc} now with $\alpha_{\tau}^{m-1}$, which gives 
all the estimates in Step 2 of the proof of
Proposition~\ref{prop-ED-Ch5-existence}. For the estimate
\eq{EUL-L-heter-Gronwall+}, it is important that, due to \eq{ass-difuse1},
$\mathscr{L}_\lambda^{}$ does not depend on $\alpha$ because we do not have
$\nabla\alpha$ estimated at this stage. For the estimate of $\alpha$
and $\mu$, we use \eq{EUL-L-viscoelastodyn+3disc-dif} tested by $\mu_{\tau}^k$,
which
yields the discrete energy balance like \eq{EUL-L-basic-engr-balance-disc-dam}.
Specifically, using the convexity of $\varphi_\lambda(X,F,\cdot)$,
testing \eq{EUL-L-viscoelastodyn+3disc-dif} by $\mu_{\tau}^k$ yields
\begin{align}\nonumber
\int_\varOmega\frac{\varphi_\lambda(\bm\xi_{\tau}^k,\Fek,\alpha_{\tau}^k)-\varphi_\lambda(\bm\xi_{\tau}^k,\Fek,\alpha_{\tau}^{k-1})}\tau
&+M(\bm\xi_{\tau}^k,\alpha_{\tau}^k)|\nabla\mu_{\tau}^k|^2\,\d\xx
\\[-.6em]&
\le-\int_\varOmega[\varphi_\lambda]_\alpha'(\bm\xi_{\tau}^k,\Fek,\alpha_{\tau}^k)\vv_{\tau}^k\Cdot\nabla\alpha_{\tau}^k\,\d\xx\,.
\label{ass-diffuse-term1}\end{align}
Furthermore, using the Green formula and the boundary condition
$\nn\Cdot\vv_{\tau}^k=0$, we obtain the estimate
\begin{align}\nonumber
-&\int_\varOmega[\varphi_\lambda]_\alpha'(\bm\xi_{\tau}^k,\Fek,\alpha_{\tau}^k)
\vv_{\tau}^k\Cdot\nabla\alpha_{\tau}^k\,\d\xx
\\[-.6em]&\nonumber\
=\int_\varOmega[\varphi_\lambda]_X'(\bm\xi_{\tau}^k,\Fek,\alpha_{\tau}^k)\Cdot
(\vv_{\tau}^k\Cdot\nabla)\bm\xi_{\tau}^k
+[\varphi_\lambda]_F'(\bm\xi_{\tau}^k,\Fek,\alpha_{\tau}^k)\Colon
(\vv_{\tau}^k\Cdot\nabla)\Fek\,\d\xx
\\[-.4em]&\nonumber\quad
\le C\big\|\vv_{\tau}^k\big\|_{L^2(\varOmega;\R^3)}\Big(
\sup\big|[\varphi_\lambda]_X'\big|
\big\|\nabla\bm\xi_{\tau}^k\big\|_{L^{r_1}(\varOmega;\R^{3})}\!\!
+\sup\big|[\varphi_\lambda]_F'\big|
\big\|\nabla\Fek\big\|_{L^r(\varOmega;\R^{3\times3})}\Big)=:C_{k,\tau}^{(1)}\,;
\end{align}
here we used \eq{est-of-xi-disc} and \eq{est-of-E-disc}.
Estimating also
\begin{align}\nonumber
\int_\varOmega&\frac{\varphi_\lambda(\bm\xi_{\tau}^{k-1},\Fekk,\alpha_{\tau}^{k-1})
-\varphi_\lambda(\bm\xi_{\tau}^k,\Fek,\alpha_{\tau}^{k-1})}{\tau\ }\,\d\xx
\\[-.6em]&\hspace{2.5em}\le
\sup\big|[\varphi_\lambda]_X'\big|
\Big\|\frac{\bm\xi_{\etau}^k{-}\bm\xi_{\etau}^{k-1}\!\!}\tau\:\Big\|_{L^1(\varOmega;\R^3)}\!\!
+\sup\big|[\varphi_\lambda]_F'\big|
\Big\|\frac{\Fek{-}\Fekk\!\!}\tau\:\Big\|_{L^1(\varOmega;\R^{3\times3})}\!\!=:C_{k,\tau}^{(2)},
\label{ass-diffuse-term3}\end{align}
summing \eq{ass-diffuse-term1}--\eq{ass-diffuse-term3} for $k=1,...,T/\tau$, we
obtain the estimate 
\begin{align}\nonumber
\!\!\!\!&\int_\varOmega\!\varphi_\lambda(\bm\xi_{\tau}^k,\Fek,\alpha_{\tau}^k)\,\d\xx
+\tau\!\sum_{m=1}^k\int_\varOmega\!M(\bm\xi_{\tau}^m,\alpha_{\tau}^m)|\nabla\mu_{\tau}^m|^2\,\d\xx
\\[-1.3em]&\hspace{15em}\le\!\int_\varOmega\!\varphi_\lambda(\bm\xi_0,\Fezero,\alpha_0)\,\d\xx
+\tau\!\sum_{m=1}^kC_{k,\tau}^{(1)}\!+C_{k,\tau}^{(2)}\,.
\end{align}
From this and from the previous estimates, we obtain also an estimate of
$\sqrt{M(\overlinexitau,\overlinealphaetau)}\nabla\overlinemuetau$ in
$L^2(I{\times}\varOmega;\R^3)$. To obtain an estimate of
$\nabla\overlinealphaetau$ in $L^2(I{\times}\varOmega;\R^3)$ from the
strategy \eq{est-of-nabla-alpha+}, one needs \eq{ass-difuse2}. Moreover,
by comparison from the former relation in
\eq{EUL-L-viscoelastodyn+3disc-dif}, we obtain the bound for
$\pdt{}\alpha_\tau={\rm div}(M(\overlinexitau,\overlinealphaetau)\nabla\overlinemutau)
-\overlinevvtau\Cdot\nabla\overlinealphaetau$ in
$L^2(I;H^1(\varOmega)^*)$. By the (generalized) Aubin-Lions theorem, we thus
obtain the strong convergence (in terms of subsequences) of $\overlinealphaetau$
in $L^2(I{\times}\varOmega)$.

The latter relation in \eq{EUL-L-viscoelastodyn+3disc-dif} means, in other
words, that $0\le\overlinealphaetau\le1$ is to satisfy the variational
inequality $(\overlinemutau-
[\varphi_\lambda]_\alpha'(\overlinexitau,\overlineFetau,\overlinealphaetau))(\wt\alpha-\overlinealphaetau)\ge0$ for any $0\le\wt\alpha\le1$, i.e.
\begin{align}\label{diffusion-VI}
\int_0^T\!\!\!\int_\varOmega\big(\overlinemutau-[\varphi_\lambda]_\alpha'(\overlinexitau,\overlineFetau,\overlinealphaetau)\big)(\wt\alpha-\overlinealphaetau)\,\d\xx\d t\ge0\,.
\end{align}
The limit passage in this variational inequality is easy when combining the
weak convergence of $\overlinemutau$ in $L^2(I{\times}\varOmega)$ and the
mentioned strong convergence of $\overlinealphaetau$.

The weak formulation of \eq{EUL-L-diff} with the initial condition
$\alpha|_{t=0}=\alpha_0$ and the mentioned boundary conditions  
leads to the integral identity
\begin{align}
&\int_0^T\!\!\!\int_\varOmega\!\alpha\pdt{\widetilde\alpha}
+\alpha({\rm div}\,\vv)\widetilde\alpha
+\big(M(\bm\xi,\alpha)\nabla\mu{+}\alpha\vv\big)\Cdot\nabla\widetilde\alpha\,\d\xx\d t
=\!\int_\varOmega\!\alpha_0\widetilde\alpha(0)\,\d\xx
\label{def-ED-Ch5-diffusion}
\end{align}
which is to hold for any $\widetilde\alpha$ smooth with $\widetilde\alpha(T)=0$
and with the inclusion $\mu\in\varphi_\alpha'(\bm\xi,\Fe,\alpha)+N_{[0,1]}(\alpha)$
to hold a.e.\ in $I{\times}\varOmega$. This needs in particular
$M(\bm\xi,\alpha)\nabla\mu\in L^1(I{\times}\varOmega;\R^3)$,
which is granted even in $L^2(I{\times}\varOmega;\R^3)$
by the boundedness of $\sqrt{M(\bm\xi,\alpha)}\nabla\mu$ in
$L^2(I{\times}\varOmega;\R^3)$ and the boundedness of $\sqrt{M(\bm\xi,\alpha)}$
in $L^\infty(I{\times}\varOmega)$ when assuming $M$ continuous on
$\barOmega{\times}[0,1]$. The above sketched arguments augmenting the proof of
Proposition~\ref{prop-ED-Ch5-existence} leads to:

\begin{proposition}\label{prop-difusion}
Let $\varphi\in C^1(\barOmega{\times}\GL_3^+{\times}[0,1])$ be bounded from below
with $\varphi(\bm\xi,\Fe,\cdot)$ convex
and let  (\ref{EUL-L-ass}b-e) with $\zetap$ modified as $\alpha$-dependent hold
for $p>r>3$ together with \eq{ass-difuse} hold with $M\in C(\barOmega{\times}[0,1])$.
Then the discrete scheme (\ref{EUL-L-viscoelastodyn+disc-dam}a-c,e) with
\eq{EUL-L-viscoelastodyn+3disc-dif}
gives approximate solutions which converge in terms of subsequences and,
for $\lambda$ sufficiently large, any such a limit is a 
weak solution of the problem (\ref{Euler-hypoplast-dam}a--d,f) and
\eq{EUL-L-diff} with the boundary conditions
\eq{Euler-small-BC-hyper} and $\nn\Cdot\nabla\mu=0$ on $\varGamma$ and
the initial conditions \eq{IC-large} together with the condition
$\alpha|_{t=0}^{}=\alpha_0\in H^1(\varOmega)$ with $0\le\alpha_0\le1$ on $\varOmega$.
\end{proposition}

For an analysis of such a model (but  without inelastic strain and 
discretized by the Galerkin method), we refer to \cite{RouSte23VESS}.

\subsection{ Diffusion of the extensive variable (a diffusant content)}\label{sec-diff-ext}

An alternative model  to Sect.~\ref{sec-diffuse}  could consider 
$\alpha$ as an {\it extensive variable} (i.e., a diffusant content instead
of concentration) valued in $[0,+\infty)$ and governed by equation
\begin{align}\label{EUL-L-diff+}
\DT\alpha={\rm div}\big(M(\bm\xi,\alpha)\nabla\mu)-({\rm div}\,\vv)\alpha
\ \ \text{ with }\ \mu=\varphi_\alpha'(\bm\xi,\Fe,\alpha)
\end{align}
instead of \eq{EUL-L-diff} while the Cauchy stress  $\TT$ 
in the momentum equation \eq{Euler1-hypoplast-dam} would be augmented by the
pressure-like term 
$\alpha\varphi_\alpha'(\bm\xi,\Fe,\alpha)\bbI$  which compensates
the term arising from $({\rm div}\,\vv)\alpha$ while testing
\eq{EUL-L-diff+} by $\mu$ to the energy balance \eq{EUL-L-energy-balance-diff},
cf.\ also \cite{MieRou25GTMF} and \cite{PesRou26CPVE}. The energetics is again
as \eq{EUL-L-energy-balance-diff}, from which the a-priori estimates for
$(\varrho,\vv,\pp,\Fe)$ and $\mu$ follows. In addition, adopting again the
boundary condition $\nn\Cdot\nabla\mu=0$, 
the total diffusant content
$\int_\varOmega\alpha(t)\,\d\xx=\int_\varOmega\alpha_0\,\d\xx$ is constant so
that $\alpha\in L^\infty(I;L^1(\varOmega))$.
Unlike in Section~\ref{sec-diffuse},
the inequality $\alpha\ge0$ is now guaranteed without any indicator function
in \eq{EUL-L-diff+} by assuming the degeneracy $M(\bm\xi,0)=0$. 
From this, we can also see the uniform bound both for
$\|\Fe(t)\|_{L^\infty(\varOmega;\R^{3\times3})}$ and also for
$\|\varphi(\bm\xi(t),\Fe(t),\alpha(t))\|_{L^1(\varOmega)}$, which is then used
for the choice of the truncation threshold $\lambda<+\infty$.

Not to create more coupling than in the scheme in Sect.~\ref{sec-diffuse}, this
additional pressure must be discretized and truncated carefully - here
``nonlocally''. Together with (\ref{EUL-L-viscoelastodyn+disc-dam}a,c,e),
we devise
\begin{subequations}\label{EUL-L-viscoelastodyn+disc-dam++}
\begin{align}
&\nonumber
\!\!\frac{\pp_\etau^k{-}\pp_\etau^{k-1}\!\!}\tau\,=
{\rm div}\Big(\mathscr{T}_\lambda(\bm\xi_{\etau}^{k-1},\Fekk,\alpha_{\eetau}^{k-1})
+\DD_\etau^k-\pp_\etau^k{\otimes}\vv_\etau^k
\\[-.3em]&
\hspace*{8.5em}
+\frac{\alpha_{\tau}^{k-1}[\varphi_\lambda]_\alpha'(\bm\xi_{\tau}^{k-1},\Fekk,\alpha_{\tau}^{k-1})}{1+(\|\varphi_\lambda(\bm\xi_{\tau}^{k-1},\Fekk,\alpha_{\tau}^{k-1})\|_{L^1(\varOmega)}^{}-\lambda)^+}
\bbI\Big)
+\varrho_\etau^k\GRAVITY_{\tau}^k\ \ \text{ and}
\label{EUL-L-viscoelastodyn+2disc-dif+}
\\[-.2em]
\label{EUL-L-viscoelastodyn+3disc-dif+}
&\!\frac{\!\alpha_{\tau}^k{-}\alpha_{\tau}^{k-1}\!\!\!}\tau
={\rm div}\big(M(\bm\xi_{\tau}^k,\alpha_{\tau}^k)\nabla\mu_{\tau}^k
-\vv_\tau^k\Cdot\nabla\alpha_{\tau}^k\big)
\ \,\text{ with }\
\mu_{\tau}^k\!=[\varphi_\lambda]_\alpha'(\bm\xi_{\tau}^k,\Fek,\alpha_{\tau}^k)\,.\!
\end{align}
\end{subequations}
For the test \eq{EUL-L-viscoelastodyn+2disc-dif+} by $\vv_\etau^k$ while using
\eq{EUL-L-viscoelastodyn+0disc-dam}, the additional truncated stress
needs the assumption
\begin{subequations}\label{ass-difuse+}
\begin{align}\label{ass-difuse1+}
&
C_\lambda:=\sup_{X\in\varOmega,\ F\in\R^{3\times3},\ \alpha\ge0}\!\frac{
\alpha\big|[\varphi_\lambda]_\alpha'(X,F,\alpha)\big|}{
1+\varphi_\lambda(X,F,\alpha)}<\infty
\intertext{
so that we can estimate this truncated pressure term tested by $\vv_\tau^k$ as}
\nonumber
\int_\varOmega&{\rm div}\Big(
\frac{\alpha_{\tau}^{k-1}[\varphi_\lambda]_\alpha'(\bm\xi_{\tau}^{k-1},\Fekk,\alpha_{\tau}^{k-1})}{1+(\|\varphi_\lambda(\bm\xi_{\tau}^{k-1},\Fekk,\alpha_{\tau}^{k-1})\|_{L^1(\varOmega)}^{}-\lambda)^+}
\bbI\Big)\,\vv_\tau^k\,\d\xx
\\&\nonumber=-
\frac{\int_\varOmega\alpha_{\tau}^{k-1}[\varphi_\lambda]_\alpha'(\bm\xi_{\tau}^{k-1},\Fekk,\alpha_{\tau}^{k-1})\,{\rm div}\,\vv_\tau^k\,\d\xx}{1+(\|\varphi_\lambda(\bm\xi_{\tau}^{k-1},\Fekk,\alpha_{\tau}^{k-1})\|_{L^1(\varOmega)}^{}-\lambda)^+}
\\&\nonumber\le
\frac{\|\alpha_{\tau}^{k-1}[\varphi_\lambda]_\alpha'(\bm\xi_{\tau}^{k-1},\Fekk,\alpha_{\tau}^{k-1})\|_{L^1(\varOmega)}\|{\rm div}\,\vv_\tau^k\|_{L^\infty(\varOmega)}}{1+(\|\varphi_\lambda(\bm\xi_{\tau}^{k-1},\Fekk,\alpha_{\tau}^{k-1})\|_{L^1(\varOmega)}^{}-\lambda)^+}
\le\frac23 C_\lambda(1{+}\lambda)\|{\rm div}\,\vv_\tau^k\|_{L^\infty(\varOmega)}\,.
\intertext{Thus we obtain all the
estimates in Step~2 of the proof of Proposition~\ref{prop-ED-Ch5-existence}.
The a-priori estimate for $\alpha$ is now obtained slightly differently than
as in Section~\ref{sec-diffuse} by executing the test 
\eq{EUL-L-viscoelastodyn+3disc-dif+} by $\alpha_\etau^k$, for which we need a
modification of the assumption \eq{ass-difuse2} as}
&\sup_{X\in\varOmega,\ F\in\R^{3\times3},\ \alpha\ge0}\!\!\!
M(X,\alpha)\big|[\varphi_\lambda]_{(X,F),\alpha}''(X,F,\alpha)\big|
+\frac1{M(X,\alpha)[\varphi_\lambda]_{\alpha\alpha}''(X,F,\alpha)}
<\infty\,,
\label{ass-difuse2+}\end{align}
\end{subequations}
Relying on the convexity of the function $\alpha\mapsto\alpha^2$, this test
results to the inequality 
\begin{align}\nonumber
&\int_\varOmega\!\frac{\alpha_\etau^k}2\,\d\xx
+\!\sum_{m=1}^k\int_\varOmega\! M(\bm\xi_{\tau}^m,\alpha_{\tau}^m)
[\varphi_\lambda]_{\alpha\alpha}''(\bm\xi_{\tau}^m,\Fem,\alpha_{\tau}^m)|\nabla\alpha_{\tau}^m|^2\,\d\xx
\\[-.3em]&\nonumber\
\le\!\int_\varOmega\!\frac{\alpha_0}2\,\d\xx
-\!\sum_{m=1}^k\int_\varOmega\!\alpha_{\tau}^m\vv_{\tau}^m\Cdot\nabla\alpha_{\tau}^m
\\[-.2em]&\qquad\ \ \nonumber
+
M(\bm\xi_{\tau}^m,\alpha_{\tau}^m)\Big([\varphi_\lambda]_{\bm\xi\alpha}''(\bm\xi_{\tau}^m,\Fem,\alpha_{\tau}^m)\Cdot\nabla\bm\xi_{\tau}^m
+[\varphi_\lambda]_{\Fe\alpha}''(\bm\xi_{\tau}^m,\Fem,\alpha_{\tau}^m)\Colon\nabla\Fem\Big)\Cdot\nabla\alpha_{\tau}^m
\,\d\xx
\\[-.2em]&\nonumber\
\le\frac12\|\alpha_0\|_{L^2(\varOmega)}^2+\!\!\sum_{m=1}^k\bigg(\int_\varOmega\!\frac{{\rm div}\,\vv_{\tau}^m}2|\nabla\alpha_{\tau}^m|^2\,\d\xx+C\|\nabla\bm\xi_{\tau}^m\|_{L^r(\varOmega;\R^{3\times3})}^{}\!
+C\|\nabla\Fem\|_{L^r(\varOmega;\R^{3\times3\times3})}^{}\bigg)
\end{align}
with $C$ related with \eq{ass-difuse2+}.
From this, by the discrete Gronwall inequality, we can obtain the bound of
$\overlinealphaetau$ and $\underlinealphaetau$
in $L^\infty(I;L^2(\varOmega))\,\cap\,L^2(I;H^1(\varOmega))$
for sufficiently small time steps $\tau>0$. These bounds on $\nabla\overlinealphaetau$
and $\nabla\underlinealphaetau$ together with the resulting
bound on $\pdt{}\alpha_\tau={\rm div}(M(\overlinexitau,\overlinealphaetau)
\nabla\overlinemutau-\overlinevvtau\Cdot\nabla\overlinealphaetau)$ with
$\overlinemutau=
[\varphi_\lambda]_\alpha'(\overlinexitau,\overlineFetau,\overlinealphaetau)$
in $L^1(I;W^{1,1}(\varOmega)^*)$ causes the strong convergence (in terms of
subsequences) of $\overlinealphaetau$ and $\underlinealphaetau$ for $\tau\to0$
by the Aubin-Lions theorem. Here we used boundedness of
$M(\overlinexitau,\overlinealphaetau)
\nabla\overlinemutau$ in $L^1(I{\times}\varOmega)$ for which we needed
$M(\overlinexitau,\overlinealphaetau)$ bounded in $L^1(I{\times}\varOmega)$
for which we need the growth condition $M(X,\alpha)\le C(1+\alpha^6)$.

The weak formulation of \eq{EUL-L-diff+} with the initial condition
$\alpha|_{t=0}=\alpha_0$ and the mentioned boundary conditions  
leads to the integral identity
\begin{align}
&\int_0^T\!\!\!\int_\varOmega\!\alpha\pdt{\widetilde\alpha}
+\big(M(\bm\xi,\alpha)\nabla\mu{+}\alpha\vv\big)\Cdot\nabla\widetilde\alpha\,\d\xx\d t
=\!\int_\varOmega\!\alpha_0\widetilde\alpha(0)\,\d\xx
\ \ \ \text{ with }\ \mu=\varphi_\alpha'(\bm\xi,\Fe,\alpha)
\label{def-ED-Ch5-diffusion+}
\end{align}
holds for any $\widetilde\alpha$ smooth with $\widetilde\alpha(T)=0$.
The above sketched arguments augmenting the proof of
Proposition~\ref{prop-ED-Ch5-existence} leads to a modification of
Proposition~\ref{prop-difusion} as:

\begin{proposition}
Let $\varphi\in C^1(\barOmega{\times}\GL_3^+{\times}\R^+)$ be bounded from below
with $\varphi(\bm\xi,\Fe,\cdot)$ convex
and let  (\ref{EUL-L-ass}b-e) with $\zetap$ modified as $\alpha$-dependent hold
for $p>r>3$ together with \eq{ass-difuse1} and \eq{ass-difuse+} hold with
$M\in C(\barOmega{\times}[0,+\infty))$ non-negative with $M(\bm\xi,0)=0$ and
with the growth condition $M(X,\alpha)\le C(1+\alpha^6)$.
Then the discrete scheme
(\ref{EUL-L-viscoelastodyn+disc-dam}a,c,e) with \eq{EUL-L-viscoelastodyn+disc-dam++}
gives approximate solutions which converge in terms of subsequences and,
for $\lambda$ sufficiently large, any such a limit is a  
weak solutions to the problem (\ref{Euler-hypoplast-dam}a--d,f) and
\eq{EUL-L-diff+} with the boundary conditions
\eq{Euler-small-BC-hyper} and $\nn\Cdot\nabla\mu=0$ on $\varGamma$ and
the initial conditions \eq{IC-large} together with 
the condition $\alpha|_{t=0}^{}=\alpha_0\in H^1(\varOmega)$ with $\alpha_0\ge0$ on $\varOmega$.
\end{proposition}

The condition \eq{ass-difuse2+} well complies with the degeneracy
$M(\bm\xi,0)=0$ which guarantees non-negativity of $\alpha$. The canonical
choice is $M(\bm\xi,\alpha)=m_0(\bm\xi)\alpha$ and then $\varphi$ contains a
convex term $\phi(\alpha):=\alpha({\rm ln}(\alpha/\alpha_{\rm eq})-1)$
with a given equilibrium concentration $\alpha_{\rm eq}>0$ at which
$\phi(\cdot)$ attains minimum, so that
$\varphi_{\alpha\alpha}''$ leads to the term $\phi''(\alpha)=1/\alpha$ which
exactly compensates the degeneracy in the mobility $M$ and thus
\eq{ass-difuse2+} can be satisfied. Also it complies with \eq{ass-difuse1+}
since $\alpha\phi'(\alpha)=\alpha{\rm ln}(\alpha/\alpha_{\rm eq})=\phi(\alpha)+\alpha$ is dominated by (the corresponding contribution to) the energy,
namely $\phi(\alpha)+\alpha\le
(1{+}\alpha_{\rm eq})\phi(\alpha)+\alpha_{\rm eq}^2{\rm e}^{1/\alpha_{\rm eq}}$.

{\small
\section*{{\large Acknowledgments}}

\vspace*{-1em}
This research has been partially supported also from the CSF (Czech Science
Foundation) project 23-06220S
and by the institutional support RVO: 61388998 (\v CR).

\medskip
}

\baselineskip10pt

{\small

}

\end{document}